\documentclass[11pt]{article}
\usepackage{natbib}
\usepackage{amsmath, amsfonts, amsthm, amssymb, mathtools}
\usepackage{authblk, color, bm, bbm, graphicx, epstopdf}
\usepackage[labelfont = bf, labelsep = period]{caption} 
\usepackage{subcaption, xspace, enumitem,adjustbox}
\usepackage{booktabs, tabularx, multirow, rotating, array,lscape,hyperref,makecell}
\usepackage{longtable,tabu}
\usepackage[bottom,hang]{footmisc}
\usepackage[flushleft]{threeparttable}
\usepackage{accents}

\makeatletter
\newcommand{\mylabel}[2]{#2\def\@currentlabel{#2}\label{#1}}
\makeatother
\usepackage{tikz,pgfplots}
\usetikzlibrary{positioning,arrows.meta}
\usetikzlibrary{shapes.geometric, arrows,calc}
 \tikzstyle{methodBox} = [rectangle, rounded corners, minimum width=1cm, minimum height=1.75cm,text centered, draw=black, very thick, text width=5cm]
  \tikzstyle{arrow} = [thick,->,>=stealth, line width=2pt]
\usepackage{pgfplots}
\usetikzlibrary{intersections}
\usetikzlibrary{shapes, arrows, chains}
\usetikzlibrary{arrows.meta}

\newcommand{\x}{\boldsymbol{x}}

\newtheorem{obs}{Observation}

\usepackage{fullpage}[2cm]
\title{Constructing Tight Quadratic Relaxations for Global Optimization: II. Underestimating Difference-of-Convex (D.C.) Functions}

\author{William~R.~Strahl$^{1,2}$}
\author{Arvind~U.~Raghunathan$^3$}
\author{Nikolaos~V.~Sahinidis$^{2,4}$}
\author{Chrysanthos~E.~Gounaris$^{1,2}$\thanks{Corresponding author: gounaris@cmu.edu}}

\affil{\small $^1$Department of Chemical Engineering, Carnegie Mellon University, Pittsburgh, PA, 15213, USA \\ 
			  $^2$Center for Advanced Process Decision-making, Carnegie Mellon University, Pittsburgh, PA, 15213, USA \\ 
			  $^3$Mitsubishi Electric Research Labs, Cambridge, MA, 02139, USA \\ 
			  $^4$H.~Milton Stewart School of Industrial~\& Systems Engineering and School of Chemical~\& Biomolecular Engineering, Georgia Institute of Technology, Atlanta, GA, 30332, USA}

\date{}

\begin{document}

\maketitle

\begin{abstract}
Recent advances in the efficiency and robustness of algorithms solving convex quadratically constrained quadratic programming (QCQP) problems motivate developing techniques for creating convex quadratic relaxations that, although more expensive to compute, provide tighter bounds than their classical linear counterparts. In the first part of this two-paper series~\citep{quadue1}, we developed a cutting plane algorithm to construct convex quadratic underestimators for twice-differentiable convex functions, which we extend here to address the case of non-convex difference-of-convex (d.c.) functions as well. Furthermore, we generalize our approach to consider a hierarchy of quadratic forms, thereby allowing the construction of even tighter underestimators. On a set of d.c. functions extracted from benchmark libraries, we demonstrate noteworthy reduction in the hypervolume between our quadratic underestimators and linear ones constructed at the same points. Additionally, we construct convex QCQP relaxations at the root node of a spatial branch-and-bound tree for a set of systematically created d.c. optimization problems in up to four dimensions, and we show that our relaxations reduce the gap between the lower bound computed by the state-of-the-art global optimization solver BARON and the optimal solution by an excess of 90\%, on average. 

\end{abstract}

\section{Introduction}\label{sec:introduction}
Global optimization is utilized in a variety of research fields, including recent applications in route planning for unmanned air vehicles~\citep{ozturk2023biobjective}, energy storage system selection, design, and operation~\citep{zantye2023theseus}, heat exchanger network synthesis~\citep{zhou2024novel}, resource recovery from wastewater~\citep{durkin2024surrogate}, pathway optimization using kinetic metabolic models for mammalian cells~\citep{lu2023kinetic}, profit maximization of scheduling hydrogen production with solar power and grid energy supply~\citep{yang2022scheduling}, and minimizing nitrogen oxide emissions produced by incinerating explosive waste materials~\citep{kim2022data}, to name but a few. The interested reader can find a plethora of additional historical applications of global optimization in eleven different research areas nicely organized in Table~2 of~\citet{boukouvala2016global}. Indeed, the ubiquitous utilization of global optimization in the literature manifests its substantial impact as a tool for scientific discovery and also highlights the importance of any advancements in this area.

The key feature of global optimization algorithms is their guarantee for identification of a global--rather than local--optimal solution. In contrast to optimization problems that are convex (i.e, the objective function and the feasible set are convex), for which a local minimum is also global~\citep{bertsekas2003convex}, non-convex problems exhibit additional complexity in the sense that non-convex objective surfaces and constraints and/or disjoint feasible regions lead to the existence of multiple extrema. For such problems, many global optimization algorithms, such as spatial branch-and-bound, compute rigorous lower and upper bounds on the objective value and successively refine those on partitioned space until the bounds converge (within some $\varepsilon$~tolerance) to a global optimal solution. Without question, the quality of the bounds impacts the convergence of these algorithms.

For minimization problems, the objective function evaluated at any feasible solution provides an upper bound, while the globally optimal objective value of a relaxation of the problem provides a lower bound. Convex relaxations find extensive applications in these algorithms for determining lower bounds because their local (and hence global) minimum are efficient to compute; also, the quality of the convex relaxation used (i.e., how closely the convex relaxation approximates the original problem) directly impacts the quality of the computed bound. To this end, an enormous amount of research has focused on developing tight convex relaxations for various specific problem structures, which has entailed significant work on deriving convex envelopes for many types of functions. Additionally, after the required expense to construct the relaxations, algorithms available to solve the convex relaxation, such as linear programming (LP) or non-linear programming (NLP) solvers, contribute to the effort required to determine the bound, and thereby contribute greatly to the overall efficiency of the algorithm. Noting the improvement in the robustness and efficiency of algorithms that solve convex quadratically constrained quadratic programs (QCQPs)~\citep{hansmittlemann}, we focus in this work on creating tight convex quadratic relaxations of non-convex optimization problems. In particular, we extend our previous methodology for constructing quadratic outer approximations of twice-differentiable convex functions~\citep{quadue1} to the case of non-convex difference-of-convex (d.c.) functions, which is a very general class of functions arising ubiquitously in global optimization applications.

In the literature, quadratics are frequently utilized to create convex relaxations of non-convex functions, but in most cases they do not result in an actual \textit{quadratic relaxation}. The $\alpha$BB methodology, for example, uses quadratics to create convex relaxations for general non-convex functions by adding a sufficiently large convex quadratic term to overcome the non-convexities of the function over the entire domain~\citep{maranas1995finding,androulakis1995alphabb,adjiman1998global,adjiman1998global}. The approach has been generalized in~\citet{akrotirianakis2004role,skjal2012generalization} to include modifications to the diagonal and nondiagonal terms of the quadratic, extended in~\citet{meyer2005convex} to use subintervals of the domain to create a spline from piecewise quadratic functions, and utilized in~\citet{gounaris2008tight,gounaris2008tightb} to create tight piecewise linear convex relaxations for non-convex functions. While the $\alpha$BB methodology utilizes quadratics in the construction of underestimators, it only creates a quadratic underestimator if the function that is underestimated is itself quadratic. Excluding this special case, none of the $\alpha$BB variants discussed above produce create convex quadratic relaxations. 

Other studies in the literature have proposed quadratic underestimators for convex functions, but these works restrict construction of their underestimators to specific classes of functions or problems. \citet{su2018improved} proposed scaling the second-order term of a Taylor series approximation at a point of construction and provided a procedure for determining the tightest scaling parameter for underestimating a restrictive class of functions, which has been the foundational study for our own work. \citet{buchheim2013quadratic} use quadratic underestimators in integer programming, introducing quadratic cuts that are generated for a specific class of functions where a matrix, $Q$, is known a~priori to satisfy $\nabla^2 f(x) \succcurlyeq Q$ for all $x \in \mathbb{R}^n$. \citet{olama2023distributed} includes quadratic cuts in a mixed integer nonlinear programming algorithm derived specifically for strongly convex functions where the strong convexity parameter is easily determined.
Last, there exist studies in the literature that construct quadratic underestimators for general non-convex functions, but these procedures either lack guarantees of convexity for the quadratic underestimators, or they only apply in restricted contexts. In particular, \citet{streeter2022automatically,streeter2023sharp} introduce quadratic over- and underestimators for general functions based on so-called Taylor polynomial enclosures, which are efficiently constructed even for functions of high dimension, but are not guaranteed to be convex. In their work, \citet{ouanes2015new} present a procedure to create a convex quadratic underestimator by subtracting a quadratic perturbation from a linear interpolant of a general non-convex function, but their method only applies in the restricted context of minimizing a single multivariate function subject to box constraints. 

As discussed, in this paper we will focus on underestimating d.c. functions. The importance of such functions in optimization is evidenced by their extensive appearance in the literature, which is due to their general applicability, succinctly captured by Tuy's statement that \textit{``every continuous global optimization problem on a compact set can be reformulated as a d.c. optimization problem''} in~\citet{horst2013handbook} (pp.~149-150). In fact, early pioneers established that any twice-differentiable continuous function defined over a convex set in $\mathbb{R}^n$ is representable as a d.c. function~\citep{hartman1959functions}, and a later study showed that any piecewise linear continuous function can be expressed as a d.c. function too~\citep{melzer1986expressibility}. Difference-of-convex functions naturally arise in a diverse number of specific applications, including but not limited to problems in signal processing, communications, and networking~\citep{gasso2009recovering,alvarado2014new,nguyen2023dca}, transportation~\citep{holmberg1999production}, facility location~\citep{chen1998solution}, image processing~\citep{lou2015weighted}, clustering~\citep{bagirov2016nonsmooth,bagirov2018nonsmooth}, and machine learning~\citep{le2017dca,awasthi2024dc,askarizadeh2023difference}. 

Relaxations for d.c. functions exist in the literature, but they are not convex quadratic, except in limited cases. These relaxations include methodologies that apply to general non-convex functions, such as $\alpha$BB~\citep{adjiman1998global,adjiman1998globalb}, McCormick relaxations~\citep{mccormick1976computability,mitsos2009mccormick,tsoukalas2014multivariate}, and the reformulation-linearization technique~\citep{sherali2013reformulation}, to name but a few. Recently, the reformulation-perspectification technique published in~\citet{bertsimas2023novel} also creates convex relaxations for d.c. functions; however, none of the aforementioned techniques construct convex quadratic relaxations for general d.c. functions, which distinguishes our work from prior relevant studies. 

In summary, we propose in this manuscript a methodology to construct convex quadratic underestimators for twice-differentiable d.c. functions. Our contributions to the literature include:
\begin{itemize}
	\item An extension of the cutting plane algorithm in~\citet{quadue1} that yields convex quadratic relaxations for d.c. functions.
	\item A hierarchy of quadratic forms established to generate relaxations of increasing tightness (at a computational cost).
	\item A computational experiment demonstrating, both qualitatively and quantitatively, the constructed underestimators on a set of functions extracted from optimization benchmark libraries.
	\item A comparison of relaxations for d.c. optimization problems at the root node of a spatial branch-and-bound tree, constructed using our methodology versus the state-of-the-art, which showcases the quality of our relaxations.
\end{itemize}

The remainder of the paper is organized as follows. In Section~\ref{sec:Methodology}, we present our proposed methodology, including (i) the extension of the algorithm in~\citet{quadue1} to construct convex quadratic underestimators for non-convex d.c. functions, and (ii) the hierarchy established by generalizing the quadratic form and introducing a shift term, as necessary. In Section~\ref{sec:computational_study}, we present our computational experiments, including (i) qualitative evidence and quantitative results for underestimators generated by our methodology, and (ii) a comparison of convex quadratic relaxations for systematically created d.c. optimization problems against relaxations created by the \texttt{BARON} solver~\citep{sahinidis1996baron}. Finally, in Section~\ref{sec:conclusions}, we offer some conclusions and state possible directions for future work.

\section{Methodology} \label{sec:Methodology}
Henceforth, we follow a notation convention where we denote vectors in boldface and matrices in capital letters. Given a d.c. function $f : \mathbb{R}^n \mapsto \mathbb{R}$, defined over a box domain $\mathcal{B} := \{\x \in \mathbb{R}^n : x^L_i \le x_i \le x^U_i \, \forall \, i = 1,\ldots,n\}$, we define a quadratic underestimator of $f$ as $q(\x; \alpha, \x_0) : \mathbb{R}^n \mapsto \mathbb{R}$, where $\x_0 \in \mathcal{B}$ is a chosen \textit{point of construction} and $\alpha$ is a scaling parameter required to realize underestimation over the full domain $\x \in \mathcal{B}$.

\subsection{Synopsis of Prior Work}\label{sec:synop}
In the first part of this two-paper series~\citep{quadue1}, we constructed a cutting plane algorithm that determined the tightest value of the scaling parameter $\alpha$ (a \textit{scalar}) of the quadratic underestimator (\ref{eq:quadform}), introduced in~\citep{su2018improved}, for general twice-differentiable convex functions. The cutting plane algorithm, as presented in~\citet{quadue1}, solves the partial epigraph reformulation of (\ref{eq:opt1}), where the objective value provides a bound on the magnitude that the quadratic overestimates the function at any point in the domain, given a value for $\alpha$. 
\begin{equation}
\label{eq:quadform}
q(\x; \alpha, \x_0) := f(\x_0) + \nabla f(\x_0)(\x-\x_0) + \frac{1}{2}(\x-\x_0)^\top \alpha\nabla^2 f(\x_0)(\x-\x_0) 
\end{equation}
\begin{equation}
\label{eq:opt1}
 \min\limits_{\x \in \mathcal{B}} f(\x) - q(\x; \alpha, \x_0)
\end{equation}

Additionally, we proposed a natural extension of the cutting plane algorithm for underestimating nonlinear functions in the context of optimization problems that feature linear constraints, which we exploited to tighten the generated underestimators by allowing the latter to overestimate in infeasible regions.
In this work, we introduce modifications at two specific parts of the methodology: (i) in the treatment of formulation (\ref{eq:opt1}), and (ii) in Step~5 of the cutting plane algorithm, where the quadratic underestimator is checked for possible overestimation at enumerated vertices and corrective action is taken, if necessary.

\subsection{Extension to Non-convex D.C. Functions}\label{sec:dc_ext}
We extend the cutting plane algorithm, originally devised for twice-differentiable convex functions, to accommodate non-convex d.c. functions. In particular, we observe that (\ref{eq:opt1}) is a d.c. optimization problem, and by defining $f(\x)$ as a d.c. function, i.e., $f(\x) := h(\x) - g(\x)$, where $h$ and $g$ are both convex, (\ref{eq:opt1}) becomes (\ref{eq:opt2}), where $g(\x) + q(\x; \alpha, \x_0)$ is convex.
\begin{equation}
\label{eq:opt2}
 \min\limits_{\x \in \mathcal{B}} h(\x) - \big[g(\x) + q(\x; \alpha, \x_0)\big]
\end{equation}

Consequently, we can employ a partial epigraph reformulation to produce (\ref{eq:opt3}), for which the objective $t - \big(g(\x) + q(\x; \alpha, \x_0)\big)$ is a concave function defined over a convex feasible set, and the same cutting plane algorithm in~\citet{quadue1} executes using only a modified objective. 
\begin{equation}
\label{eq:opt3}
\begin{array}{cl}
 \min\limits_{\x \in \mathcal{B}, t \in \mathbb{R}} & t - \big[g(\x) + q(\x; \alpha, \x_0)\big] \\ 
 \text{s.t.}             & h(\x) - t \le 0
\end{array}
\end{equation}

We make the following remarks: (i) the globally optimal objective value of (\ref{eq:opt3}) for a valid underestimator is $0$, and is attained at $\x = \x_0$, assuming that an underestimating quadratic of form (\ref{eq:quadform}) exists at $\x_0$; and (ii) the set of vertices maintained by the cutting plane algorithm is an outer approximation of only $h(\x)$, and thus $h(\x)$ should be the only function utilized to initialize and enumerate vertices. We note that the reformulation preserves the monotonically increasing property of lower bounds generated by the cutting plane algorithm, and thus the convergence proof in~\citet{quadue1}, which relies on the proof from~\citet{hoffman1981method}, holds for (\ref{eq:opt3}) without any alterations. Finally, we highlight that we can directly utilize this reformulation in conjunction with the extension presented in~\citet{quadue1} to exploit the possible presence of linear constraints in an optimization problem so as to construct even tighter relaxations.

\subsection{Point of Construction Selection}\label{sec:poc}
In the case of twice-differentiable convex functions, $\nabla^2 f(\x) \succcurlyeq 0$ holds for all $\x \in \text{dom}(f)$. However, in the extension introduced in Section~\ref{sec:dc_ext}, we now consider underestimating non-convex functions, where in many cases there exist $\x \in \text{dom}(f)$ such that $\nabla^2 f(\x) \not\succcurlyeq 0 $. The quadratic underestimator (\ref{eq:quadform}) directly scales $\nabla^2 f(\x_0)$ as well as, consequently, the eigenvalues of $\nabla^2 f(\x_0)$. Thus, that underestimator will not produce a convex quadratic underestimator for a non-convex function, unless $\nabla^2 f(\x_0) \succcurlyeq 0$ as a necessary (yet not sufficient) condition.

\begin{figure}[htb]
	\centering
	\includegraphics[width=.6\textwidth]{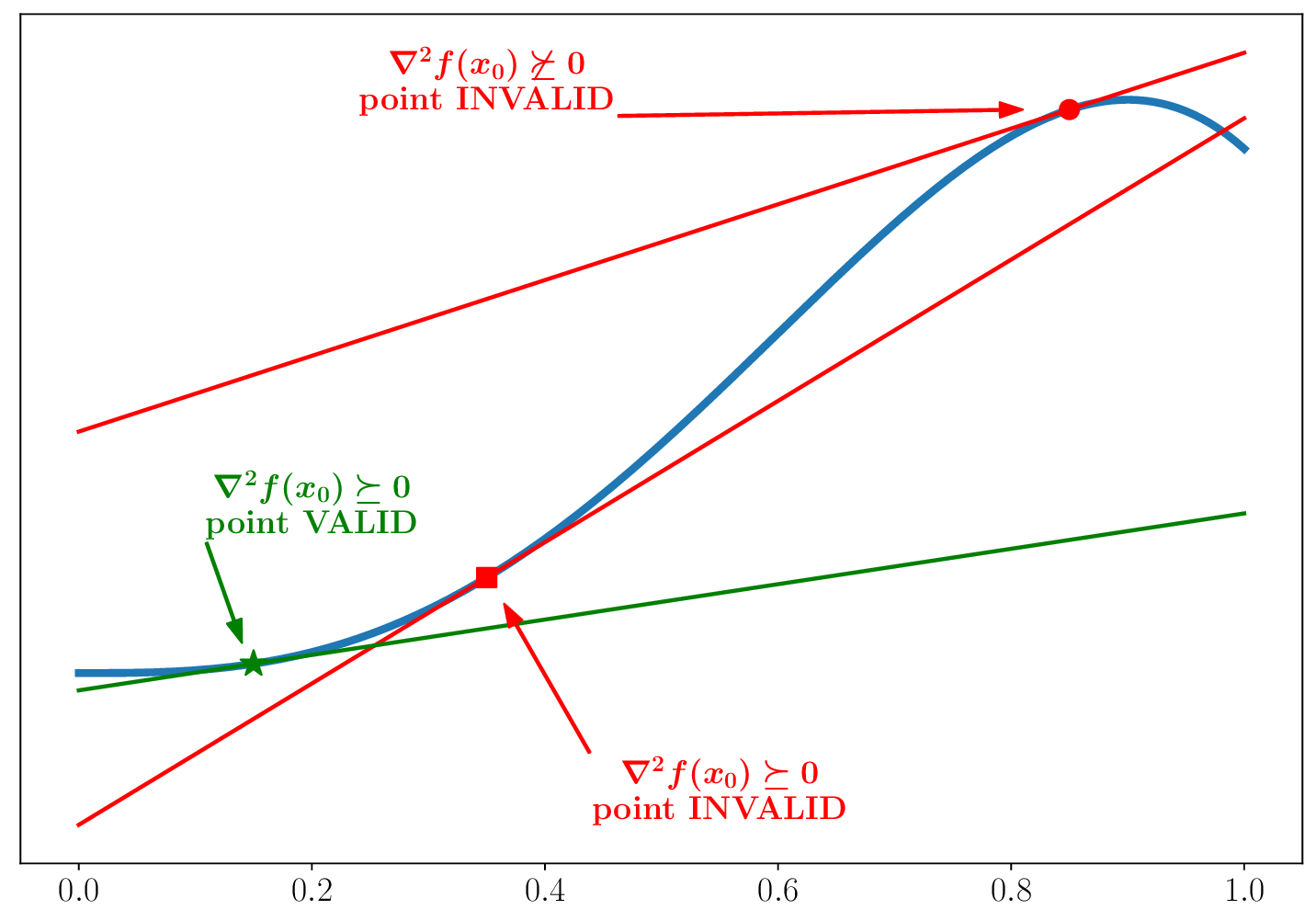}
	\caption{The univariate non-convex d.c. function $f(x) := 3x^3 - 2.5x^4$ over the domain $x \in [0,1]$ with first-order Taylor series approximations constructed at (i) $x = 0.15$, (ii) $x = 0.35$, and (iii) $x = 0.85$, illustrating respectively three different cases: (i) underestimation over the entire domain, (ii) overestimation in part of the domain despite local convexity at the point of construction, and (iii) overestimation where the point of construction is not locally convex.}
	\label{fig:POCIllustration}
\end{figure}

We also note that the validity of the bound on overestimation computed by (\ref{eq:opt3}) is predicated on the premise that $g(\x) + q(\x; \alpha, \x_0)$ is convex, and if the selection of $\x_0$ nullifies this premise, then the bound yielded by (\ref{eq:opt3}) is no longer valid. We highlight that the condition $\nabla^2 f(\x_0) \succcurlyeq 0$) is not sufficient in the non-convex case, as a first-order Taylor series approximation constructed at a point $\x_0$ (even where $\nabla^2 f(\x_0) \succcurlyeq 0$) is not guaranteed to underestimate a non-convex function over the entire domain. Figure~\ref{fig:POCIllustration} illustrates this concept with a non-convex d.c. function and three first-order Taylor series approximations, created at three candidate points of construction that each demonstrates different possibilities: (i) a point where $\nabla^2 f(\x_0) \not\succcurlyeq 0$, for which the tangent line does not even underestimate locally, (ii) a point where $\nabla^2 f(\x_0) \succcurlyeq 0$, for which the tangent line underestimates locally but not throughout the whole domain, and (iii) a point where $\nabla^2 f(\x_0) \succcurlyeq 0$ for which the tangent line constitutes a valid underestimator. Obviously, it is only in the last case that one can construct a valid quadratic underestimator of the form of (\ref{eq:quadform}).

\subsection{Generalization to a Hierarchy of Quadratic Underestimators}\label{sec:generalizationHierarchy}
Considering the fact that the scalar $\alpha$ uniformly scales the elements of $\nabla^2 f(\x_0)$ in (\ref{eq:quadform}), we generalize the quadratic form and present a hierarchy of quadratic underestimators that allows for the construction of tighter underestimators. Before presenting the methodology for each underestimator, we restate key ideas and relevant algorithmic steps for understanding and implementing our hierarchy. In the remainder of the manuscript, we use $\x^\ast$ to indicate a point at which $f(\x^\ast) - q(\x^\ast; \alpha^{(k)}, \x_0) < -\varepsilon$, where $\varepsilon$ is a user-specified tolerance and superscript ``$(k)$'' denotes the value of the scaling variable at iteration $k$ of the cutting plane algorithm.

Corrective action is taken in Step~5 of the cutting plane algorithm to suitably recompute the parameter of the quadratic. Here, we restate Step~5 as well as Observation~\ref{obs:mono} from~\citet{quadue1}, which elucidates a key property for our quadratic underestimators, i.e., that the quadratic monotonically decreases with respect to a decrease in $\alpha$, which can be exploited to improve algorithmic efficiency.

\par{Step~5.} \textit{Evaluate Underestimation:}
\begin{center} $\alpha^{(k)} \gets \min\limits_{(\x_v, t_v) \in \mathcal{H}^+} \left\{\frac{2\left(f(\x_v) - (f(\x_0) + \nabla f(\x_0)(\x_v-\x_0)\right)}{(\x_v-\x_0)^\top \nabla^2f(\x_0)(\x_v-\x_0)} \ : \ f(\x_v) - q(\x_v; \alpha^{(k)}, \x_0) < -\varepsilon \right\}$,
\end{center}
where $\mathcal{H}^+$ is the set of new vertices enumerated with the introduction of a cutting plane. 

\begin{obs}[for proof, see~\citet{quadue1}]
	\label{obs:mono} 
	The quadratic underestimator in (\ref{eq:quadform}) monotonically decreases as $\alpha$ decreases; that is, $q(\x; \alpha_1, \x_0) \le q(\x; \alpha_2, \x_0)$ for all $\alpha_1 < \alpha_2$.
\end{obs}

The goal of the following presentation is to generalize the quadratic underestimator form to allow the construction of tighter underestimators, while preserving the efficiency of the algorithm by enforcing the quadratic to monotonically nonincrease each time the parameters are updated. We note that changes at $\x^\ast$ are the same, regardless of whether the function is convex or d.c., and will therefore hereafter reference the function to be underestimated as $f(\x)$, which in the d.c. case can be substituted directly with $h(\x) - g(\x)$. We begin our presentation with the original method developed in~\citet{quadue1}, which for purposes of the hierarchy, we shall refer to as the ``scalar'' method.

\subsubsection{Scalar Method}\label{sec:S}
In the scalar method (also labeled method ``S''), we decrement the scalar $\alpha$ of the quadratic form (\ref{eq:quadform}) using the update rule~\ref{eq:S}. Here, we note the monotonic relationship of the quadratic with $\alpha$, as per Observation~\ref{obs:mono}, which implies that $\frac{1}{2}(\x^\ast-\x_0)^\top \nabla^2 f(\x_0)(\x^\ast-\x_0) \ne 0$ for all $\x^\ast$ such that $f(\x^\ast) - q(\x^\ast; \alpha^{(k)}, \x_0) < -\varepsilon$.
\begin{equation}
\label{eq:S}
\alpha^{(k)} \gets \frac{2\big(f(\x^\ast) - (f(\x_0) + \nabla f(\x_0)(\x^\ast-\x_0)\big)}{(\x^\ast-\x_0)^\top\nabla^2f(\x_0)(\x^\ast-\x_0)}\tag{S}
\end{equation}

\subsubsection{Diagonal Method}\label{sec:D}
We generalize the quadratic form in method ``S'' by introducing a diagonal matrix of scaling parameters, $A \in \mathbb{R}^{n \times n}$, into the eigenvalue decomposition of $\nabla^2 f(\x_0)$, $Q \Lambda Q^{-1}$, as shown in (\ref{eq:quadform2}), which permits the scaling parameters, i.e., the diagonal elements of $A$, to modify each eigenvalue of $\nabla^2 f(\x_0)$ differently, while the off-diagonal elements are fixed to 0. Hence, we refer to this method as the ``diagonal'' method, or method ``D''.
\begin{equation}
\label{eq:quadform2}
q(\x; A, \x_0) := f(\x_0) + \nabla f(\x_0)(\x-\x_0) + \frac{1}{2}(\x-\x_0)^\top Q A \Lambda Q^{-1}(\x-\x_0) 
\end{equation}

The new quadratic form (\ref{eq:quadform2}) requires a different procedure for determining the scaling parameters, such that each new selection:
(i) ensures the quadratic underestimates the function at $\x^\ast$;
(ii) monotonically decreases the quadratic over the entire domain with each consecutive update to the parameters;
and (iii) produces the tightest possible quadratic underestimator over the entire domain, when the algorithm converges.
Observing that the (diagonal) elements of $A$ participate linearly in the quadratic form, we can compute them via the linear program (\ref{eq:D}), where we explicitly enforce (i) and (ii) in the constraints.
\begin{equation}\tag{D}
  \label{eq:D}
  \begin{array}{cll} 
   \max\limits_{A} & \sum\limits_{\boldsymbol{v} \in \mathcal{S}} q(A;\boldsymbol{v}, \x_0) \\ 
   \text{s.t.}      & q(A; \x^\ast, \x_0) \le f(\x^\ast) \\
            & A_{ii} \le A_{ii}^{(k)} & \forall \, i \\
            & A_{ii} \ge 0 & \forall \, i\\
            & A_{ij} = 0 & \forall \, i,  j : \{j \ne i\},
  \end{array}
 \end{equation}
where $A_{ii}^{(k)}$ is the incumbent value for the diagonal elements of the scaling matrix (to be updated with the optimal solution of this LP).

To accomplish (iii), the above LP aims to explicitly optimize for ``tightness'' or quality of the produced underestimator.
For this, we utilize the metric (\ref{eq:metric}), defined in~\citet{quadue1}, that quantifies the tightness of the quadratic underestimator, $q(\x)$, as the fractional reduction in hypervolume from a function to it, when compared to a linear underestimator, $\ell(\x)$, generated at the same point of construction. 
\begin{equation}
	\label{eq:metric}
	M^{q(\x)}_{\ell(\x)} = \frac{\int_{\x \in \textrm{dom}(f)} \big(q(\x)-\ell\x)\big)d\x}{\int_{\x \in \textrm{dom}(f)} \big(f(\x)-\ell(\x)\big)d\x}
\end{equation}

As we consider utilizing this metric as the objective in the LPs for our methods, we make the following observation, where $\mathcal{S}$ is a set of points $\x \in \text{dom}(f)$.
\begin{obs}
	\label{obs:TightestParameters}
	In the case that $\left|\mathcal{S}\right| \to \infty$, the constraint $q(A; \boldsymbol{v}, \x_0) \le f(v)$ for all $\boldsymbol{v} \in \mathcal{S}$, in conjunction with the objective $\max\limits_{A} \sum_{\boldsymbol{v} \in \mathcal{S}} q(A; \boldsymbol{v}, \x_0)$, produces the tightest scaling parameters by construction, thereby achieving objective (iii).
\end{obs}

Since sampling an infinite amount of points at which to evaluate underestimating distances is not practical, we discretize the space by choosing $\left|\mathcal{S}\right| = 100n$ points, where $n$ is the dimension of the problem and the points are selected via a Latin hypercube approach. To embed the metric directly in the LP, we remove constant terms from (\ref{eq:metric}) to produce the objective function in (\ref{eq:D}), which maximizes the value of the quadratic at each sampled point. However, the constraints in (\ref{eq:D}) only enforce underestimation at a single point, namely $\x^\ast$. 
This fact, considered together with the monotonicity requirement imposed on the parameters, results in the algorithm making greedy initial commitments to parameter values which, although optimal in light of $\x^\ast$, may ultimately yield sub-optimal parameter values at termination.

To address this behavior, we recommend augmenting the formulation of the first LP we ever solve with the constraints $q(A; \boldsymbol{v}, \x_0) \le f(\boldsymbol{v})$ for all $\boldsymbol{v} \in \mathcal{S}$, resulting in formulation (\ref{eq:Dfirst}), which provides an update to the parameters using a holistic treatment of the space. These additional constraints are only required for the first LP executed in the cutting plane algorithm; subsequent LPs can be relaxed to contain only the single constraint $q(A; \x^\ast, \x_0) \le f(\x^\ast)$ (formulation~\ref{eq:D}), since monotonicity guarantees underestimation for these points through the remaining execution of the algorithm.
We highlight that an interesting tradeoff arises between the number of points in set $\mathcal{S}$ and the number of LPs that have to be solved throughout the algorithm. Despite yielding a larger LP (\ref{eq:Dfirst}), which however need only be solved once at the first iteration, larger $\left|\mathcal{S}\right|$ typically leads to detecting fewer points of overestimation, creating tighter underestimators while requiring fewer LPs (\ref{eq:D}) to be solved in the long run.

\begin{equation}\tag{$D^{(1)}$}
  \label{eq:Dfirst}
  \begin{array}{cll}
  \max\limits_{A} & \sum\limits_{\boldsymbol{v} \in \mathcal{S}} q(A;\boldsymbol{v}, \x_0) \\ 
  \text{s.t.}           & q(A; \x^\ast, \x_0) \leq f(\x^\ast) \\
  			            & q(A; \boldsymbol{v}, \x_0) \le f(\boldsymbol{v}) & \forall \, \boldsymbol{v} \in \mathcal{S} \\
         & A_{ii} \le A_{ii}^{(0)} & \forall \, i \\
		          & A_{ii} \ge 0 & \forall \, i\\
        &  A_{ij} = 0 & \forall \, i,  j : \{j \ne i\},
\end{array}
\end{equation}
where we initialize $A_{ii}^{(0)} = 1$ for all rows $i =1,2,\ldots,n$.

We remark the following regarding the diagonal scaling method:
(i) (\ref{eq:Dfirst}) has $\left|\mathcal{S}\right| + 2n + 1$ constraints in the first iteration, after which (\ref{eq:D}) has $2n + 1$ constraints; thus, these LPs are small in size;
(ii) we know that the value of $1$ is a valid upper bound on each element of the diagonal of $A$, $A_{ii}$, because the quadratic would overestimate locally at the point of construction if any $A_{ii} > 1$;
and (iii) ``S'' can be viewed as a restriction of method ``D'', where additional constraints are added in (\ref{eq:D}) to achieve equivalency for all $A_{ii}$; indeed, the final parameter $\alpha$ determined by method ``S'' used for all the diagonal elements of $A$ will be feasible to (\ref{eq:D}), but due to the inexactness introduced from sampling the space (finite set $\mathcal{S}$), formulation (\ref{eq:D}) can possibly produce an inferior underestimator, the mitigation of which is encouraged by (\ref{eq:Dfirst}).

\subsubsection{Matrix Method}\label{sec:M}
We further generalize the diagonal scaling method to allow modifications to the off-diagonal elements of $A$ as well, where the quadratic underestimator form is the same as (\ref{eq:quadform2}), but we remove the restriction that the off-diagonal elements of $A$ must be $0$. We shall refer to this method as the ``matrix'' method, or method ``M''.

While the quadratic form does not change from the diagonal scaling method, the introduction of off-diagonal variables in $A$ necessitates a more involved LP to preserve the required properties of the quadratic. First and foremost, we are interested in producing \textit{convex} quadratic underestimators, for which we had trivial guarantees in previous methods, as we directly scaled eigenvalues with $\nabla f(\x_0) \succcurlyeq 0$. However, modifying the off-diagonal elements gives rise to the possibility of creating non-convex quadratics even with non-negative diagonal elements. To that end, we explicitly enforce convexity during each parameter update by requiring $A$ to be diagonally dominant, which is a sufficient condition for the positive semi-definiteness of $A$. To invoke this property, we first require symmetry on $A\Lambda$ and non-negativity of all diagonal elements $A_{ii}$, which ensures that $\left[A\Lambda\right]_{ii} \ge 0$ due to $\Lambda_{ii} \ge 0$ at the selected point of construction. Then, we impose the diagonal dominance property by requiring that $\left[A\Lambda\right]_{ii} \ge \sum_{j \ne i}|\left[A\Lambda\right]_{ij}|$ for all rows $i=1,2,\ldots,n$.

Lastly, to preserve the algorithmic efficiency of the cutting plane algorithm achieved by exploiting the monotonicity property of the quadratic (Observation~\ref{obs:mono}), we utilize the fact that $A^{(k)}\Lambda \succcurlyeq A\Lambda$ implies that $\x A^{(k)}\Lambda \x \ge \x A\Lambda \x$ for all $\x \in \mathbb{R}^n$, where $A^{(k)}$ are the incumbent values of the parameters at the time the LP is executed.
This realizes the desired monotonicity property, since the orthogonality of $Q$ in the eigenvalue decomposition preserves the eigenvalues of $A\Lambda$. For the constraint requiring $A^{(k)} \Lambda \succcurlyeq A\Lambda$, we similarly impose the diagonal dominance property on matrices $(A^{(k)}- A)\Lambda$.
Based on the above, the computation of $A$ can be achieved via formulation (\ref{eq:MC}), which can be trivially reformulated to the LP in (\ref{eq:M}) after the introduction of auxiliary variables $S\in \mathbb{R}^{n \times n}$ and $T\in \mathbb{R}^{n \times n}$.
\begin{equation} \tag{M}
	\label{eq:MC}
	\begin{array}{cll}
	\max\limits_{A} & \sum\limits_{\boldsymbol{v} \in \mathcal{S}}q(A; \boldsymbol{v}, \x_0) \\ 
	\text{s.t.}				 & q(A; \x^\ast, \x_0) \le f(\x^\ast) \\
							 & A_{ii} \le A_{ii}^{(k)} & \forall \, i \\
	 						 & A_{ij}\Lambda_{jj} = A_{ji}\Lambda_{ii} & \forall \, i, j : \{j \ne i\} \\
	 						 & A_{ii} \ge 0 & \forall \, i \\
							 & A_{ii}\Lambda_{ii} \ge \sum\limits_{j \neq i}\left|A_{ij}\Lambda_{jj}\right| & \forall \, i\\ 
							 & (A^{(k)}_{ii}-A_{ii})\Lambda_{ii} \ge \sum\limits_{j \neq i}\left|(A^{(k)}_{ij}-A_{ij})\Lambda_{jj}\right| & \forall \, i,
	\end{array}
\end{equation}

\begin{equation}\tag{M$^\prime$}
 \label{eq:M}
	\begin{array}{cll}
	\max\limits_{A, S, T} & \sum\limits_{\boldsymbol{v} \in \mathcal{S}}q(A; \boldsymbol{v}, x_0) \\ 
	\text{s.t.}				 & q(A; \x^\ast, \x_0) \le f(\x^\ast) \\
						 & A_{ii} \le A_{ii}^{(k)} & \forall \, i \\
						 & A_{ij}\Lambda_{jj} = A_{ji}\Lambda_{ii} & \forall \, i, j : \{j \ne i\} \\
	  						 & A_{ii} \ge 0 & \forall \, i\\						 
		& A_{ii}\Lambda_{ii} \ge \sum\limits_{j \neq i}S_{ij} & \forall \, i\\
							 & -S_{ij} \le A_{ij}\Lambda_{jj} \le +S_{ij} & \forall \, i, j : \{j \ne i\} \\
							 & (A^{(k)}_{ii}-A_{ii})\Lambda_{ii} \ge \sum\limits_{j \neq i}T_{ij} & \forall \, i\\
							 & -T_{ij} \le (A^{(k)}_{ij}-A_{ij})\Lambda_{jj} \le +T_{ij} & \forall \, i, j : \{j \ne i\}
	\end{array}
\end{equation}

We highlight that the LP (\ref{eq:M}) has $n^2 + 2n(n-1)$ variables and $4n + 5n(n-1) + 1$ constraints. Despite being a more involved LP than the one encountered in the diagonal method, it is still expected to be very tractable in practice. 

\subsubsection{Shift}\label{sec:shift}
In Section~\ref{sec:poc}, we established a necessary--but not sufficient--condition for the selection of point of construction to yield a valid quadratic underestimator generated by our cutting plane algorithm. Here, we introduce a term into the quadratic forms (\ref{eq:quadform}) and (\ref{eq:quadform2}) to allow for their vertical (downward) shift. We show how to properly update (\ref{eq:S}), (\ref{eq:D}) (or \ref{eq:Dfirst}) and (\ref{eq:M}) to accommodate the augmented form, and how the introduction of the shift causes the necessary condition $\nabla^2 f(\x_0) \succcurlyeq 0$ to also become sufficient for our methodology to produce valid quadratic underestimators.
For all methods, we denote the scalar variable for the shift as $\gamma \in [0,\infty)$, which is initialized to $\gamma^{(0)}=0$, and which is to be negated from the applicable quadratic form (\ref{eq:quadform}) or (\ref{eq:quadform2}) in each case.
To differentiate from the previous methods without shift, we augment the label of our methods that include the shift with the letter ``S'' (e.g., ``DS'' refers to the diagonal method \emph{with shift}).

For the scalar method, we utilize a shift only if method ``S'' computes $\alpha < 0$, whereupon we reset $\alpha \gets 0$ and update $\gamma$ throughout each iteration using (\ref{eq:SS}). While computationally inexpensive, this approach (dubbed method ``SS'') will effectively produce a linear underestimator when the shift is indeed required (i.e., a vertical shift of the first-order Taylor series approximation).
\begin{equation}
\label{eq:SS}
\gamma^{(k)} \gets q(\x^\ast, 0, \x_0) - f(\x^\ast) \tag{SS}
\end{equation}

Before we present the incorporation of shift in the methods utilizing a scaling matrix $A$ in lieu of a scalar $\alpha$, we consider a hybrid method that requires identical elements in a diagonal matrix. Consistent with the rest of the hierarchy, we will refer to this as method as ``UDS'' (i.e., ``uniform diagonal with shift'').
In this method, both scalar parameters $\alpha$ (the common element in the diagonal) and $\gamma$ can be optimized simultaneously via the LP~(\ref{eq:UDS}).
We also highlight the inclusion of the constraint $\gamma \ge \gamma^{(k)}$, which imposes the monotonicity property of the quadratic from iteration to iteration.
Importantly, this enforcement of monotonicity preserves the requisite properties for the cutting plane algorithm's convergence proof~\citep{quadue1} after the addition of the shift.
We also remark that method ``UDS'' is a generalization of both methods ``SS'' and ``S'' inasmuch as any results determined by the latter two are also feasible to ``UDS'', noting however that the relative tightness of these underestimators might be affected by the selection of the set $\mathcal{S}$.
\begin{equation}\tag{UDS}
  \label{eq:UDS}
  \begin{array}{cll}
   \max\limits_{A, \gamma} & \sum\limits_{\boldsymbol{v} \in \mathcal{S}} \big[q(A;\boldsymbol{v}, \x_0) - \gamma \big] \\ 
   \text{s.t.}             & q(A; \x^\ast, \x_0) - \gamma \le f(\x^\ast) \\
             			 & A_{ii} \le A_{ii}^{(k)} & \forall \, i\\
   					     & A_{ii} \ge 0 & \forall \, i\\
   						 & A_{11} = A_{jj} & \forall \, j : \{j \ne 1\} \\ 
            			 & \gamma \ge \gamma^{(k)},
 \end{array}
 \end{equation}
where $\gamma^{(k)}$ is the incumbent value for the shift variable (to be updated with the optimal solution of this LP).

\begin{figure}[htb]
	\centering
	\includegraphics[width=.5\textwidth]{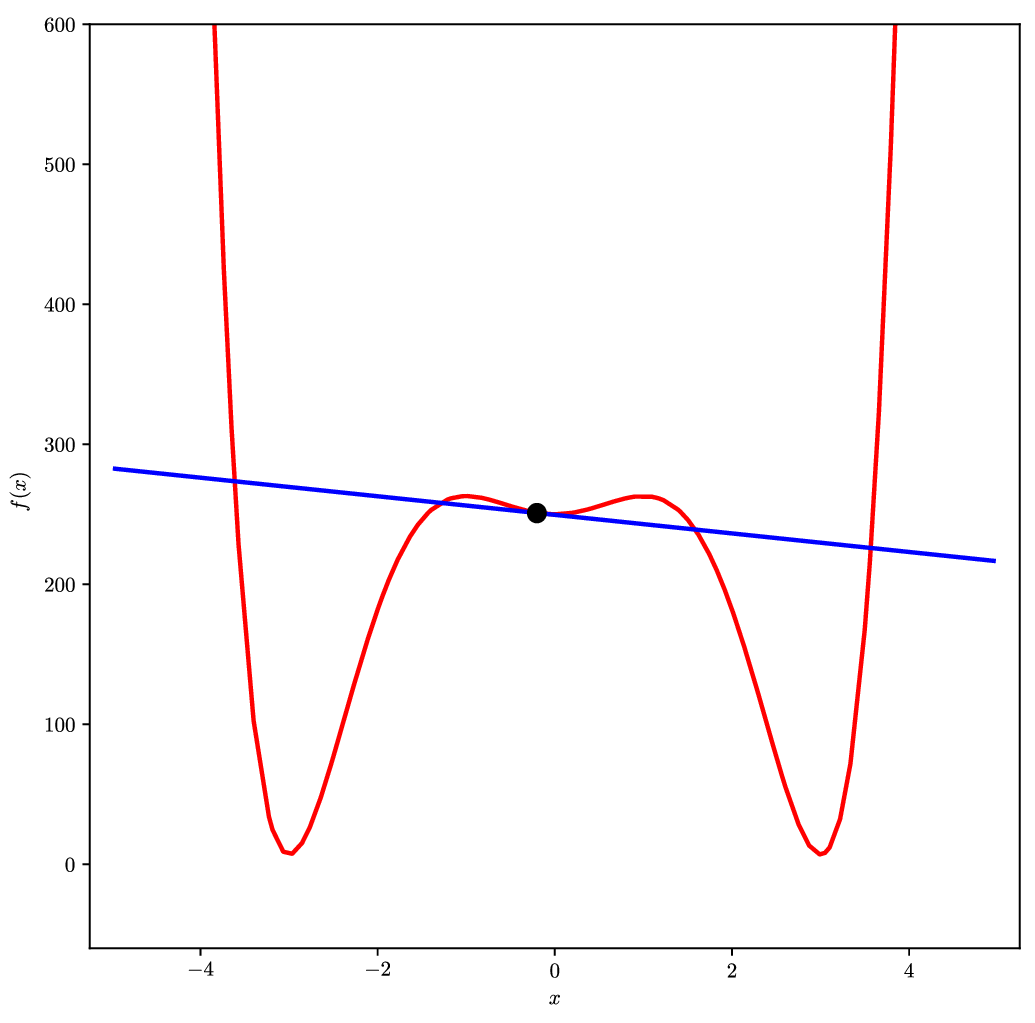}\hfill
	\includegraphics[width=.5\textwidth]{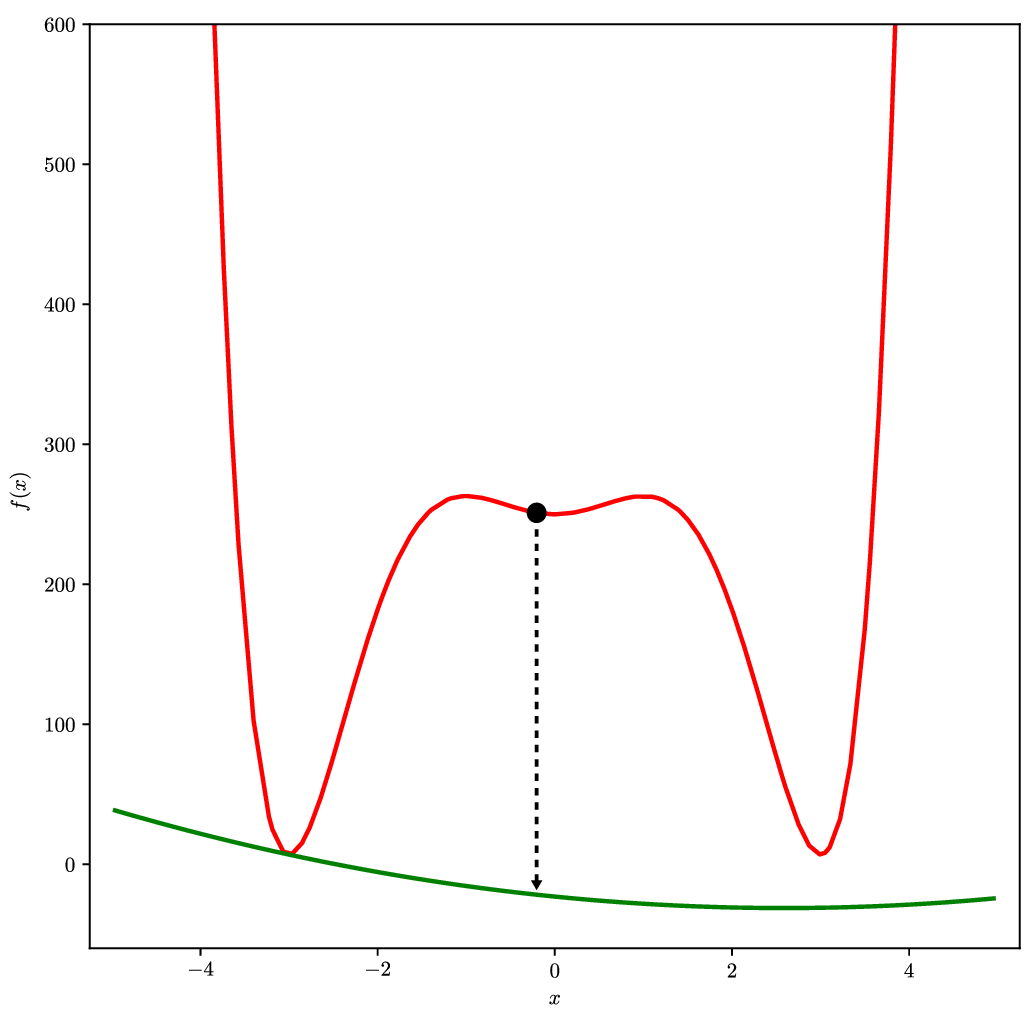}
	\caption{Introducing the shift parameter allows the construction of a valid quadratic underestimator (right, in green) for the d.c. function $f(x) = \left(27x^2 + x^6 + 250\right) - 15x^4$ at $x_0 = -0.125$ even when, at the same point of construction, a linear underestimator is not valid (left, in blue).}
	\label{fig:linearShift}
\end{figure}

In summary, in the hierarchy of using a mere scalar in the quadratic form, we have: (i) method ``S'', which may fail to produce a valid underestimator for certain points of construction in the case of non-convex d.c. functions; (ii) method ``SS'', which is efficient to compute but will only produce a linear underestimator, and should hence be reserved only for the points of construction where ``S'' fails; and (iii) method ``UDS'', which requires solving LPs but will produce a valid underestimator for all points of construction where $\nabla^2 f(\x_0) \succcurlyeq 0$ (at worst a linear underestimator, but likely a quadratic one).
Figure~\ref{fig:linearShift} demonstrates an example where utilizing method ``UDS'' constructs a successful quadratic underestimator when method ``S'' would otherwise fail due to poor choice of the point of construction.
Finally, we highlight that, in certain instances, a combination of $\gamma > 0$ and $\alpha > 1$ in method ``UDS'' could achieve even tighter underestimators, as measured by relative hypervolumes; however, for our computational studies, we keep the upper bound of parameter $\alpha$ at $1$ and initialize $\alpha^{(0)} \gets 1$.

We now turn our attention to the incorporation of shift in the methods utilizing scaling matrices, introducing methods ``DS'' and ``MS''. For this, the linear programs defined by (\ref{eq:D}) and (\ref{eq:M}) are updated analogously to (\ref{eq:UDS}), where we incorporate the variable $\gamma$ and impose the constraint to preserve monotonicity. This results in linear programs (\ref{eq:DS}) and (\ref{eq:MS}), respectively for the two methods.

\begin{equation}\tag{DS}
  \label{eq:DS}
  \begin{array}{cll}
   \max\limits_{A, \gamma} & \sum\limits_{\boldsymbol{v} \in \mathcal{S}} \big[q(A; \boldsymbol{v}, \x_0) - \gamma \big] \\ 
   \text{s.t.}             & q(A; \x^\ast, \x_0) - \gamma \le f(\x^\ast) \\
           & A_{ii} \le A_{ii}^{(k)} & \forall \, i \\
           & A_{ii} \ge 0 & \forall \, i \\
           & \gamma \ge \gamma^{(k)}
  \end{array}
 \end{equation}

\begin{equation}\tag{MS}
 \label{eq:MS}
	\begin{array}{cll}
	\max\limits_{A, \gamma, S, T} & \sum\limits_{\boldsymbol{v} \in \mathcal{S}} \big[q(A;\boldsymbol{v}, \x_0) - \gamma \big]\\ 
	\text{s.t.}				      & q(A; \x^\ast, \x_0) - \gamma \le f(\x^\ast) \\
						 & A_{ii} \le A_{ii}^{(k)} & \forall \, i \\
						 & A_{ij}\Lambda_{jj} = A_{ji}\Lambda_{ii} & \forall \, i, j : \{j \ne i\} \\
	  						 & A_{ii} \ge 0 & \forall \, i\\						 
		& A_{ii}\Lambda_{ii} \ge \sum\limits_{j \neq i}S_{ij} & \forall \, i\\
							 & -S_{ij} \le A_{ij}\Lambda_{jj} \le +S_{ij} & \forall \, i, j : \{j \ne i\} \\
							 & (A^{(k)}_{ii}-A_{ii})\Lambda_{ii} \ge \sum\limits_{j \neq i}T_{ij} & \forall \, i\\
							 & -T_{ij} \le (A^{(k)}_{ij}-A_{ij})\Lambda_{jj} \le +T_{ij} & \forall \, i, j : \{j \ne i\}\\
             				      & \gamma \ge \gamma^{(k)}
	\end{array}
\end{equation}

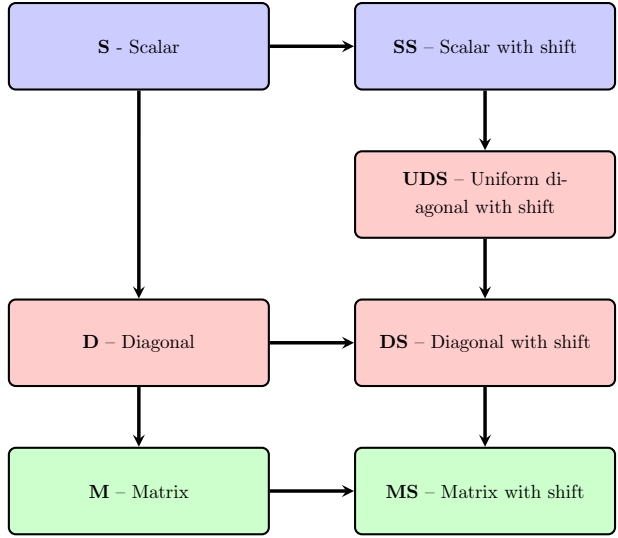
\begin{figure}[h]
   \centering
   \resizebox{0.5\textwidth}{!}{
    \begin{tikzpicture}[node distance=3cm, scale=0.2]
	\node (S) [methodBox, fill=blue!20] {\textbf{S} - Scalar};
	\node (SS) [methodBox, fill=blue!20, right of=S, xshift=4cm] {\textbf{SS} -- Scalar with shift};
	\node (DUS) [methodBox, fill=red!20, below of=SS] {\textbf{UDS} -- Uniform diagonal with shift};
	\node (DS) [methodBox, fill=red!20, below of=DUS] {\textbf{DS} -- Diagonal with shift};
	\node (MS) [methodBox, fill=green!20, below of=DS] {\textbf{MS} -- Matrix with shift};
	\node (D) [methodBox, fill=red!20, left of=DS, xshift=-4cm] {\textbf{D} -- Diagonal};
	\node (M) [methodBox, fill=green!20, left of=MS, xshift=-4cm] {\textbf{M} -- Matrix};

	\draw [arrow] (S.east) -- (SS.west);
	\draw [arrow] (S.south) -- (D.north);
	\draw [arrow] (D.south) -- (M.north);
	\draw [arrow] (SS.south) -- (DUS.north);
	\draw [arrow] (DUS.south) -- (DS.north);
	\draw [arrow] (DS.south) -- (MS.north);
	\draw [arrow] (D.east) -- (DS.west);
	\draw [arrow] (M.east) -- (MS.west);
    \end{tikzpicture}

   } 
   \caption{The hierarchy of our methods producing quadratic underestimators, based on scalar (blue), diagonal (red), and matrix (green) scalings, with the addition of the shift in those methods appearing on the right column. The arrows indicate the generation of tighter underestimators.}
	\label{fig:hierarchy}
  \end{figure}

Overall, we have presented a hierarchy of methods that can produce quadratic underestimators of non-convex d.c. functions at user-specified points of construction. This hierarchy is presented in a diagram in Figure~\ref{fig:hierarchy}, where tighter--but more computationally expensive--underestimators lie at the bottom and to the right. We remind readers that all the methods not in blue require solving an LP to update their parameters.
Finally, we emphasize the fact that all quadratic forms and update procedures are immediately amenable to, without further modification, the methodological extension presented in~\citet{quadue1}, where information of external linear constraints in optimization problems can be exploited to produce tighter underestimators by allowing them to overestimate in infeasible regions.

\section{Computational Study}\label{sec:computational_study}
We demonstrate the performance of our hierarchy of quadratic underestimators generated by our cutting plane algorithm for non-convex d.c. functions in two distinct computational experiments. In the first experiment in Section~\ref{sec:hierarchy_comparison_study}, we extract d.c. functions from optimization problems found in benchmark libraries and we construct underestimators for them, showcasing their efficiency and tightness. In the second experiment in Section~\ref{sec:root_node_relax_comp}, we compare the root node relaxation of the state-of-the-art global optimization solver \texttt{BARON} with a convex quadratic relaxation constructed via our methodology on optimization problems with d.c objectives and constraints. 

\subsection{Hierarchy Comparison Study}\label{sec:hierarchy_comparison_study}
For our study on the tightness and computational efficiency of our hierarchy of underestimators, we extracted 3 one-dimensional and 7 two-dimensional d.c. functions from the COCONUT library (\url{https://arnold-neumaier.at/glopt/coconut/Benchmark/Benchmark.html}), which compiles optimization problems from GlobalLib, CUTE, and Constraint satisfaction test problems (CSTP) libraries. 
The location of the problems in the COCONUT library and the explicit expressions used in the computational study are included in Table~\ref{tab:COCONUTDCFunctions} in the Appendix. 
We note that we employed the relationships $x_1x_2 = \frac{1}{2}\left(x_1 + x_2\right)^2 - \frac{1}{2} \left(x_1^2 + x_2^2\right)$ and $x_1^2x_2^2 = \frac{1}{2}\left(x_1^2 + x_2^2\right)^2 - \frac{1}{2}\left(x_1^4 + x_2^4\right)$ to create d.c. representations for many of the functions found in the optimization libraries. 
Where necessary, we assigned bounds to otherwise unbounded variables as noted in Table~\ref{tab:COCONUTDCFunctions}, and we scaled the range of the functions to lie within the interval $[-1,1]$ by using the scaling factor $1/\max\{|\min\limits_{\x \in \mathcal{B}}f(\x)|, |\max\limits_{\x \in \mathcal{B}} f(\x)|\}$.
For each function, we used Latin hypercube sampling of its domain until we could accumulate 25 points of construction where the function is locally convex, as assessed via eigenvalue decomposition of the Hessian at those points.
Finally, we used the value $\varepsilon = 1$e$-3$ for the convergence tolerance of the cutting plane algorithm.
Our computational experiments were executed on a machine equipped with a 1.80GHz Intel(R) Core(TM) i7-8565U~CPU running on a Ubuntu~22.04 virtual machine with 8GB~RAM and 4~logical processors, where we have implemented the algorithms in Python. To solve LPs, we used the open source HiGHS~\citep{huangfu2018parallelizing} linear program solver, which is included in the Python \texttt{SciPy} package~\citep{virtanen2020scipy}, and which was found to be adequate to address the small sized LPs encountered in this study.

\subsubsection{Results}\label{sec:hierarchy_results}
\begin{figure}[htb!]
	\centering
	\includegraphics[width=.6\textwidth]{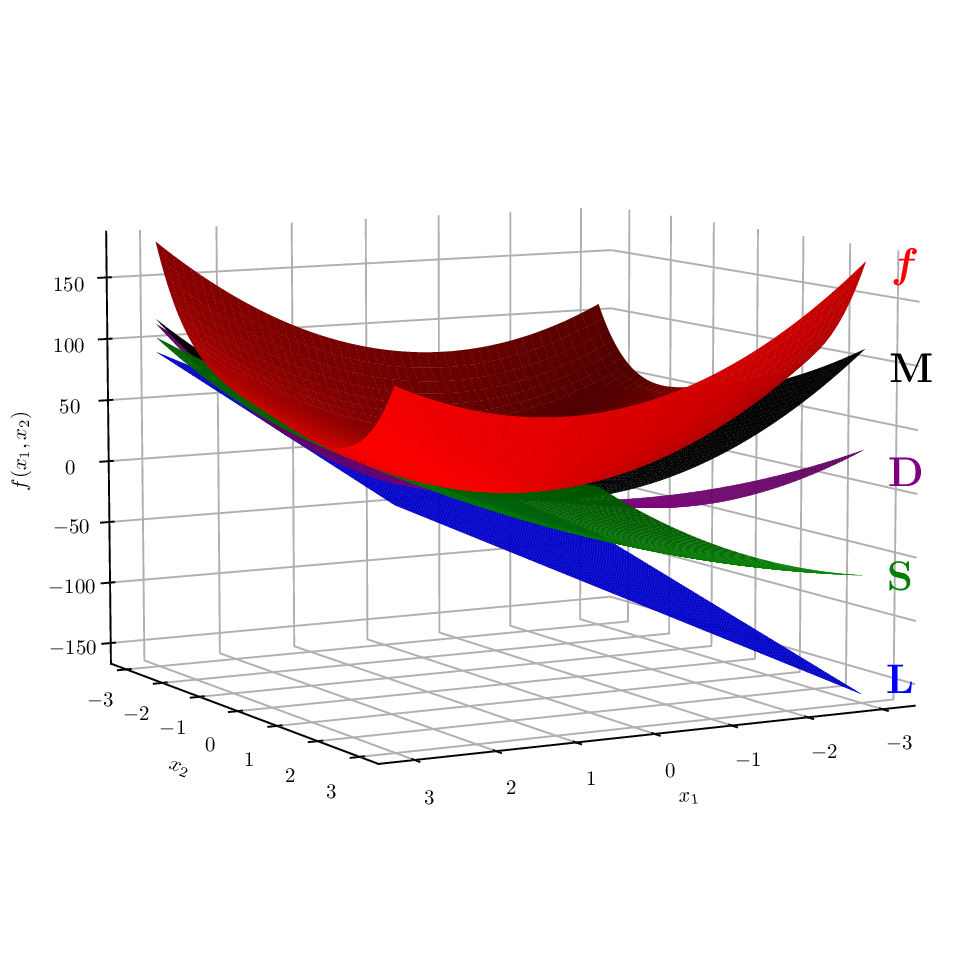}
	\caption{Qualitative evidence highlighting the improvement in the quality of underestimators by using more general methods from our hierarchy. The d.c. function $f(x_1, x_2) = \left(x_2^4 + 9x_1^2 + 2x_2^2\right) - 2\left(x_1+x_2\right)^2$ is underestimated in the box domain $[-3,3]^2$ using $(1.84, -1.04)$ as the point of construction.
	Here, ``$f$'' is the function, ``L'' is the linear underestimator, while ``S'', ``D'', and ``M'' are the quadratic underestimators constructed using the respective methods.}
	\label{fig:hierarchyPlot}
\end{figure}

\begin{figure}[htb!]
	\centering
	\includegraphics[width=.5\textwidth]{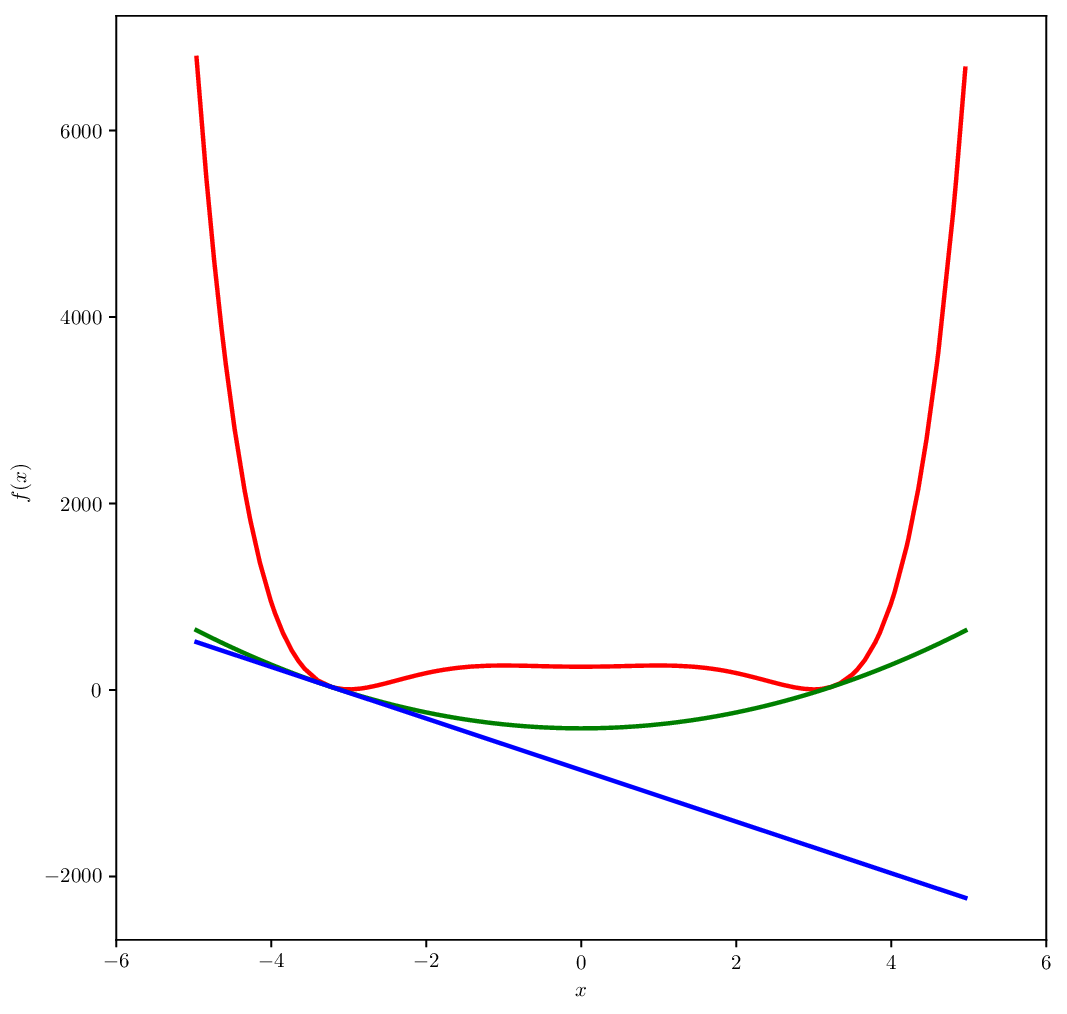}\hfill
	\includegraphics[width=.5\textwidth]{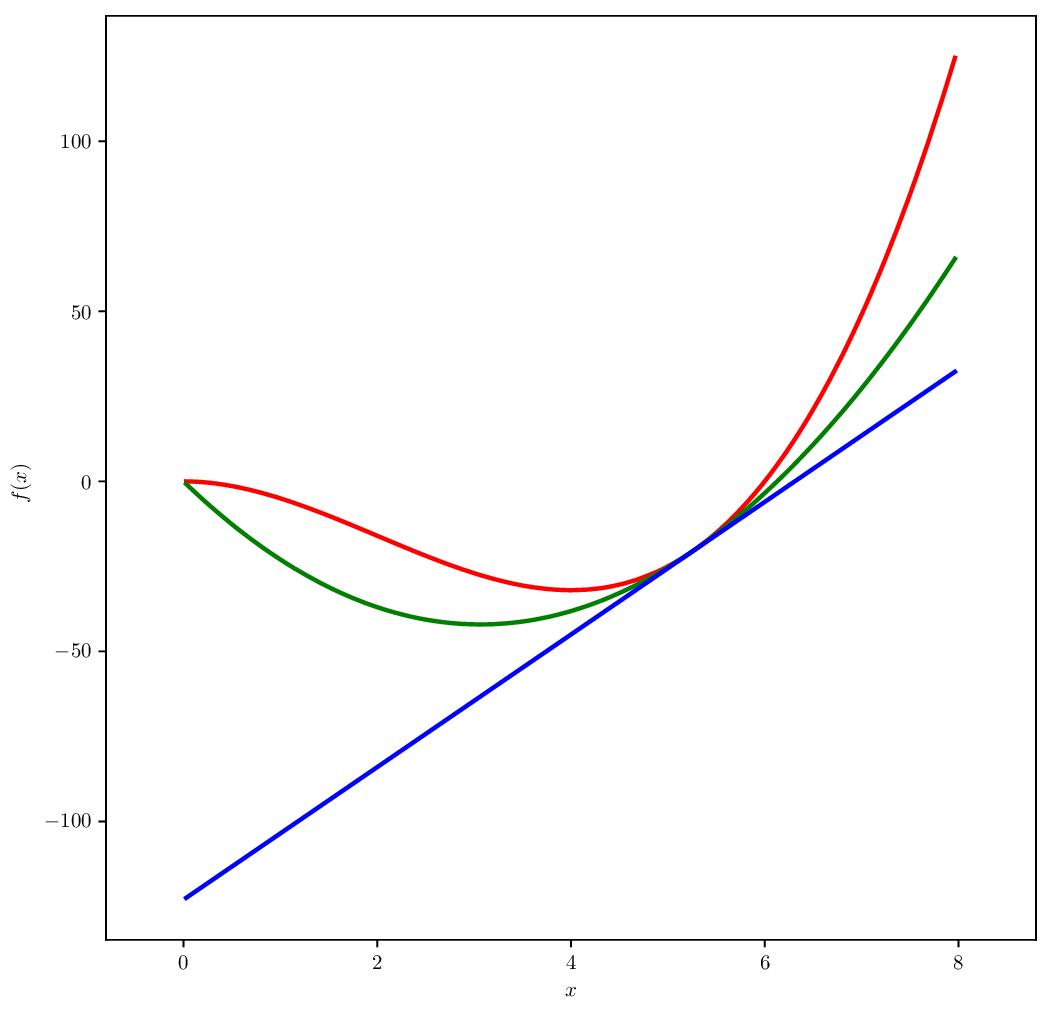}
	\caption{Quadratic underestimators (green) for the univariate d.c. functions $f(x) = \left(27x^2 + x^6 + 250\right) - 15x^4$ at $x_0 = -3.24$ (left) and $f(x) = x^3 - 6x^2$ at $x_0 = 5.24$ (right), compared to linear underestimators (blue) using the ``S'' method.}
	\label{fig:1DQuadUnderestimators}
\end{figure}

\begin{figure}[htb!]
	\centering
	\includegraphics[width=.5\textwidth]{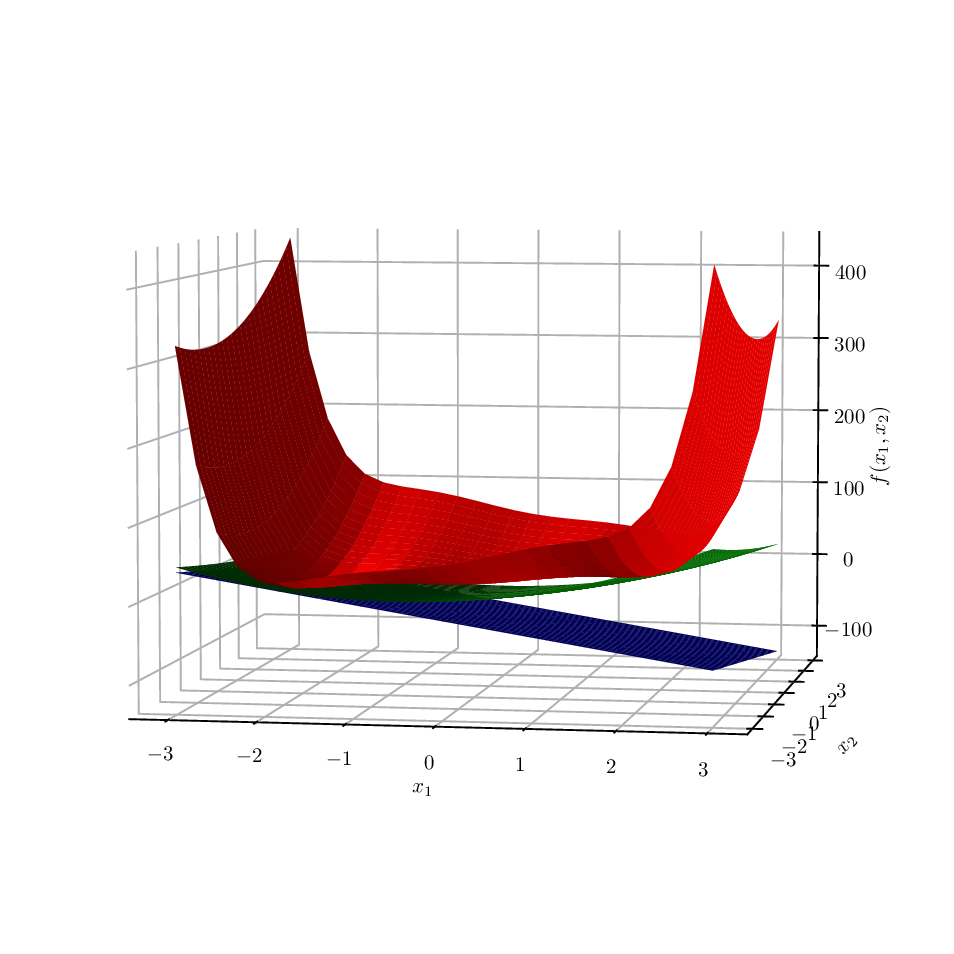}\hfill
	\includegraphics[width=.5\textwidth]{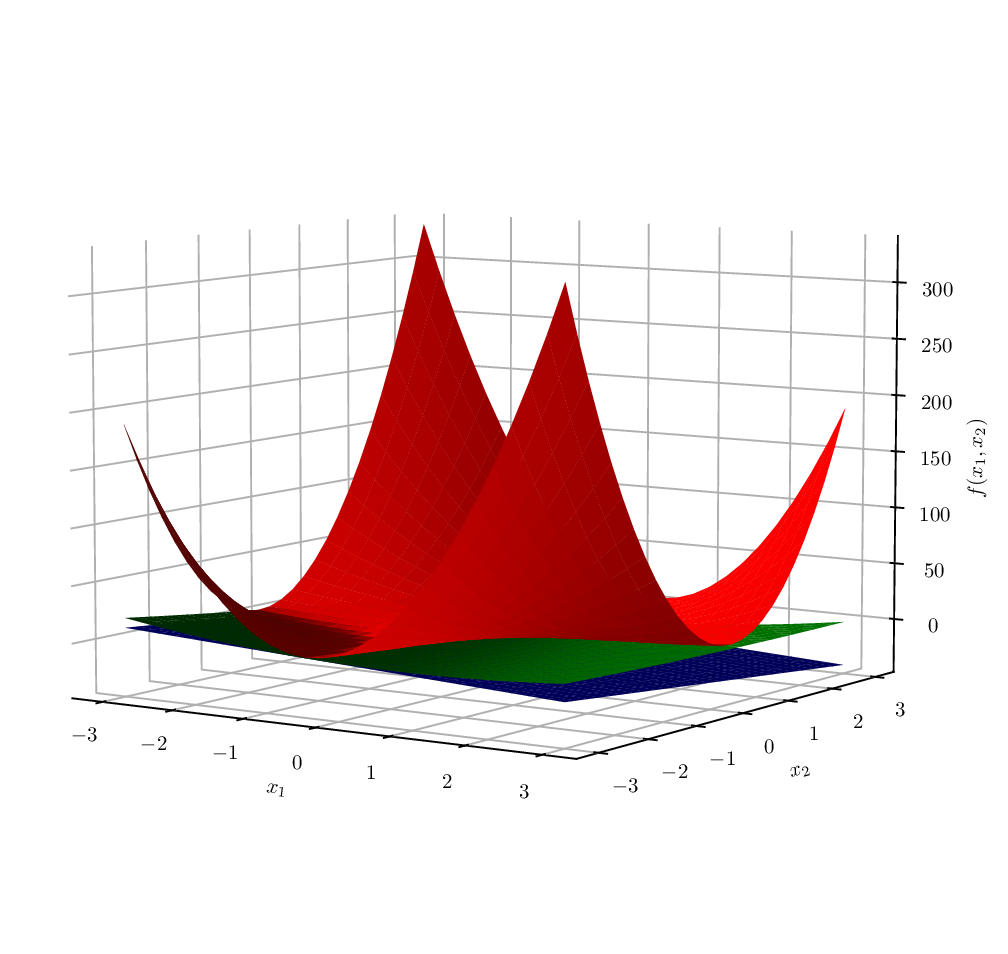}
	\caption{Quadratic underestimators (green) for the bivariate d.c. functions $f(\x) = \left(15.0x_1^2 + 9x_2^2 + x_1^6\right) - \left[3(x_1+x_2)^2 + 6.3x_1^4\right]$ at $\x_0 = (-1.95, -1.81)$ (left) and $f(\x) = \left[5x_1^2 + 5x_2^2 + \frac{3}{2}(x_1^2+x_2^2)^2\right] - \left[4(x_1 + x_2)^2 + \frac{3}{2}x_1^4 + \frac{3}{2}x_2^4\right]$ at $\x_0 = (-0.68, -2.21)$ (right), compared to linear underestimators (blue) using the ``S'' method.}
	\label{fig:2DQuadUnderestimators}
\end{figure}

Qualitatively, Figure~\ref{fig:hierarchyPlot} depicts the underestimators created by the methods ``S'', ``D'' and ``M'' for a function from the \texttt{dipigri} optimization problem, highlighting the possible improvement in underestimator tightness by using ever more sophisticated techniques.
Additionally, Figures~\ref{fig:1DQuadUnderestimators} and~\ref{fig:2DQuadUnderestimators} visualize some examples of our quadratic underestimators constructed for non-convex d.c. functions, as generated by the ``S'' method.

Performing a quantitative comparison among the methods in the hierarchy requires categorizing the functions along with points of construction into two groups: (i) those for which methods that do not utilize shift can construct a successful quadratic underestimator, where the entire hierarchy can be compared; and (ii) those for which construction procedures require the shift for success (i.e., first-order Taylor series approximations do not underestimate across the full domain of interest), and hence, where only the methods ``SS'', ``UDS'', ``DS'' and ``MS'' can be compared.
Accordingly, Tables~\ref{tab:hierarchyValidPoints} and~\ref{tab:hierarchyInvalidPoints} display the average tightness metric, computed via (\ref{eq:metric}), and the CPU time (in milliseconds) required to compute the underestimators for groups (i) and (ii), respectively.

\begin{table}[h!]
    \caption{\label{tab:hierarchyValidPoints}Computational results for different methods in the hierarchy, for points of construction that admit a valid underestimator without the shift.}
    \centering
    \resizebox{\linewidth}{!}{%
    \begin{tabular}{ c c c c c c c c c c c c c c c c c } 
     \toprule
     \multirow{2}{*}{Dimension} & \multirow{2}{*}{\# Functions} & \multirowcell{2}{\# Points of \\ Construction} & \multicolumn{7}{c}{Avg. Metric} & \multicolumn{7}{c}{Avg. CPU (ms)}\\
     \cmidrule(lr){4-10}\cmidrule(lr){11-17}
     & & & S & D & M & SS & UDS & DS & MS & S & D & M & SS & UDS & DS & MS \\ 
     \midrule
    1 & 3 & 36 & 0.559 & 0.559 & 0.559 & 0.559 & 0.559 & 0.559 & 0.559 &  4 & 11 & 13 &  6 &  11 &  11 &  14\\ 
	2 & 7 & 90 & 0.384 & 0.449 & 0.464 & 0.384 & 0.467 & 0.499 & 0.502 & 35 & 93 & 125 & 44 & 98 & 114 & 133\\	
     \bottomrule
    \end{tabular}%
    }    
\end{table}

 \begin{table}[h!]
 	\caption{\label{tab:hierarchyInvalidPoints}Computational results for different methods in the hierarchy, for points of construction that require the shift to admit a valid underestimator.}
	\centering
	\adjustbox{width=0.73\linewidth}{
	\begin{tabular}{ c c c c c c c c c c c } 
	 \toprule
	 \multirow{2}{*}{Dimension} & \multirow{2}{*}{\# Functions} & \multirowcell{2}{\#  Points of \\ Construction}& \multicolumn{4}{c}{Avg. Metric}& \multicolumn{4}{c}{Avg. CPU (ms)}\\
	 \cmidrule(lr){4-7}\cmidrule(lr){8-11}
	 & & & SS & UDS& DS& MS & SS & UDS& DS& MS\\ 
	 \midrule
	 1 & 3 & 39 & 0.000 & 0.014 & 0.014 & 0.014 &  4 &  9 &  8 & 12\\ 
	 2 & 7 & 85 & 0.000 & 0.094 & 0.115 & 0.122 & 34 & 75 & 73 & 96\\ 
	 \bottomrule
	\end{tabular}}
\end{table}

For functions of one dimension, Table~\ref{tab:hierarchyValidPoints} clearly shows that all methods perform similarly for points of construction that admit valid underestimators without requiring a shift, which is simply explained by the equivalency of the methods in the univariate case. In the case of bivariate functions, however, Table~\ref{tab:hierarchyValidPoints} shows that the shift provides some flexibility to create tighter underestimators. This is evidenced by the sharp increase of 0.083 (21.6\%) in the average metric between methods ``S'' and ``UDS'', as compared to the less significant improvement of 0.065 (16.9\%) between ``S'' and ``D''. Using the matrix methods in lieu of the diagonal ones yields only marginal improvements in underestimator tightness, as indicated by the increase of 0.015 (3.3\%) and 0.003 (0.6\%) in the average metric, respectively for the unshifted (``D'' vs. ``M'') and shifted (``DS'' vs. ``MS'') methods. As expected, method ``MS'' constructs the tightest underestimators, but it also exacts the greatest computational expense.

Table~\ref{tab:hierarchyInvalidPoints} demonstrates the effectiveness of the shift and the capability for the methodologies to construct \textit{quadratic} underestimators at points where, without the shift, the methods would fail.
First, we clearly state that all methods with the shift successfully create underestimators for all points of construction, without exception. However, given that the first-order Taylor series approximation constructed at these points is not a valid underestimator, we ought to use an altered definition of the metric to convey tightness.
More specifically, in Table~\ref{tab:hierarchyInvalidPoints}, the reported metric compares the linear underestimator generated by method ``SS'' with underestimators constructed by the other methods; hence, method ``SS'' is reported in Table~\ref{tab:hierarchyInvalidPoints} as having a metric of $0$, by definition. Notably, the other methods show a metric greater than $0$, indicating that the final underestimators generated by these methods retain some curvature. Indeed, for the bivariate case, we observe that simultaneously optimizing the parameter(s) with the shift produces underestimators that reduce the volume between the function and the linear underestimators generated by ``SS'' by 9.4\%, 11.5\%, and 12.2\%, respectively for methods ``UDS'', ``DS'' and ``MS''. These results demonstrate the importance of using the shift parameter to successfully achieve quadratic underestimation at any and all points of construction for which $\nabla^2 f(\x_0) \succcurlyeq 0$ and to improve the tightness of these underestimators.

Tables~\ref{tab:hierarchyDetailsValid} and~\ref{tab:hierarchyDetailsInvalid} in the Appendix provide detailed data for each method regarding the average and standard deviation of cutting plane algorithm iterations, number of vertices enumerated, metric, CPU time, and number of LP solves for the computational study.
Finally, we remark that, in our extraction of d.c. functions, certain pathological cases (e.g., a function from the \texttt{sisser} optimization problem showed in Figure~\ref{fig:sisserPlot} in the Appendix were identified but not removed from our benchmark set. In these cases, \textit{pathological} refers to functions for which the cutting algorithm will fail to generate a successful quadratic underestimator at all points $\{\x \in \text{dom}(f) : \nabla^2 f(\x) \succcurlyeq 0\}$ without the shift.

\subsection{Root Node Relaxation Comparison with State-of-the-art}\label{sec:root_node_relax_comp}
In this study, we demonstrate the quality of relaxations constructed using our quadratic underestimation methodology on a set of d.c. optimization problems.
In particular, we employ method ``DS'', which affords a good tradeoff between quality and tractability, to construct quadratic underestimators for each non-linear function in these problems, setting the tolerance of the cutting plane algorithm to $\varepsilon = 1$e$-3$.
We then solve the resulting convex QCQP relaxation using \texttt{IPOPT}~v3.14.13~\citep{wachter2006implementation} with the linear solver \texttt{MUMPS}~v5.6.2~\citep{MUMPS1,MUMPS2} to determine a lower bound at the root node of each problem, which we compare with the lower bound at the root node computed by \texttt{BARON}~v24.5.8~\citep{sahinidis1996baron}. We use all default settings of \texttt{BARON}, including the option \texttt{NoutPerVar = 4} that dictates the number of outer approximations per variable for convex multivariate functions. For a fair comparison, we mirror this by choosing to construct four quadratic underestimators per dimension at points of construction determined using Latin hypercube sampling, for each non-linear function present in the benchmark problems. Aside from this, we use the same settings for the cutting plane algorithm, as presented in Section~\ref{sec:hierarchy_comparison_study}.
We also highlight that we include all auxiliary techniques available by default in \texttt{BARON} (e.g., bounds tightening), even though equivalent techniques are not utilized in our methodology. Despite such an unfavorable setup, however, we will demonstrate the superior quality of the lower bounds computed using our convex quadratic relaxation for problems defined in dimensions greater than one.

\subsubsection{Generation of D.C. Optimization Problem Library}\label{sec:opt_generation}
Noting the scarcity of benchmark d.c. optimization problems available in the literature, we systematically created optimization problems that are defined in one to four dimensions and that include linear, convex, and (non-convex) d.c. functions. Whereas the explicit problem instances used in our computational study are provided in the Appendix, in this section we outline the general procedure we employed in creating these instances.

We generate optimization problems parameterized by a tuple, $\left(n, m_\ell, m_c, m_{dc}\right)$, where $n$ is the dimension of the problem, $m_\ell$ is the number of linear constraints, $m_c$ is the number of convex constraints, and $m_{dc}$ is the number of non-convex d.c. constraints. For our study, we produce instances defined using all combinations of $n \in \{1,2,3,4\}$, $m_\ell = 1$ if $n > 1$ (no linear constraints added in univariate optimization problems), $m_c \in \{1,2\}$ and $m_{dc} \in \{1,2,3\}$. In total, we produce 24 problems (6 for each dimension). In all cases, we restrict the domain of the variables to $[-1,1]$.

To create a new instance, we first assign $n$ to the desired dimension, and then select a d.c. function to serve as the objective, followed by $m_\ell$, $m_c$, and $m_{dc}$ linear, convex, and d.c. functions, respectively, to serve as constraints. 
As we add a new constraint to the problem, we ensure it is not redundant (i.e., it does reduce the feasible space) by utilizing a set of $100n$ points selected via Latin hypercube sampling to check for feasibility. More specifically, given the domain of the variables, $\x \in [-1,1]^n$, we determine a~priori the upper and lower bounds of the range of each function, where the upper bound defines a right-hand side for the constraint to encompass the full feasible space (i.e., makes the constraint redundant), and the lower bound defines a right-hand side that eliminates the entire feasible space (i.e., causes the problem to be infeasible). Via binary search, we pick the right-hand side to use for this constraint such that $20$\% of the remaining sampled feasible points are eliminated. For example, after we sample 400 points for a four dimensional problem, the first constraint added will eliminate 80 of those from the feasible space, while the next constraint will eliminate an additional 64 ($20$\% of the remaining 320) points, resulting in 256 feasible points in the original sample.

Consequently, the order in which the constraints are added to the optimization problem is important: we add the linear constraint first, followed by the convex constraints, and add the d.c. constraints last. We highlight that adding the constraints in this order maximizes the contribution of the non-convex d.c. functions in defining the feasible space, which aligns with the goal of this work. 

We generate linear constraints for use in our optimization problems by considering the linear form $\sum_{i=1}^n \beta_i x_i \le \phi$, where $\beta_i = 1$ or $\beta_i = -1$ (chosen randomly) for each dimension $i$. As per the process explained above, we identify bounds on range of each linear expression and execute a binary search on $\phi$ until the feasible space (as represented by the sample points) has been reduced by $20$\%.
In regards to generating convex functions, we do so by randomly sampling from a core set of convex functional forms extracted from benchmark libraries and expanded using established rules that preserve convexity.
To create d.c. functions, we take special care to ensure the functions are non-convex and multi-modal over the prescribed domain. We begin with a core set of six univariate non-convex d.c. functions (see Table~\ref{tab:funcUsed} in the Appendix), which we add up to construct higher dimensional d.c. functions while also adding a ``linking term'', $\mathcal{L}_i(\x) := \frac{1}{2p}(\sum_{j=1}^p x_j)^4$, where $p$ is the number variables included in a particular multi-dimensional d.c. function. The linking term induces inseparable variable dependencies on the function output, which prevents trivially separating terms that could otherwise have been underestimated separately. After randomly shuffling all possible combinations of higher dimensional d.c. functions utilizing the core set of six, we use the first d.c. function in the objective, and subsequent functions in the d.c. constraints until we have added a total of $m_{dc}$ constraints.

While we omit many specific details from our problem construction description, we include the explicit problems in Tables~\ref{tab:1D1rootNodeProblems}--\ref{tab:4D1rootNodeProblems} in the Appendix for reference.
We note that these tables only present floating point numbers to three decimals of precision; consequently, we supply the \texttt{.nl} files defining the problems with full decimal precision as executed in this study as Supplementary Material to this paper.
We additionally provide a comprehensive list of the points of construction for each function of the study, included as comments in each file. 

\subsubsection{Results}\label{sec:root_node_comparison_results}
Table~\ref{tab:rootNodeStudy} provides aggregate results for the root node relaxation comparison, where we report the number of instances for which each alternative approach--QCQP relaxation versus \texttt{BARON}--computes a superior root node lower bound as well as the average reduction of the gap between the optimal solution and the \texttt{BARON} lower bound that the QCQP relaxation affords us.
These results demonstrate that the convex QCQP relaxation produces substantially superior lower bounds compared to \texttt{BARON}, reducing the root node relaxation gap by an excess of 92\%, on average, for problems defined in two or more dimensions. For univariate problems, the QCQP relaxation produces tighter lower bounds than \texttt{BARON} at the root node for five out of the six problem instances, while for the sixth problem, \texttt{BARON} reduced the gap between the QCQP bound and the optimal solution by only 11\%.  
These results provide strong evidence of the tightness of the convex QCQP relaxation constructed by our methodology.
We acknowledge that, while this comparison provides evidence for the quality of the quadratic underestimators generated by our methodology, the latter requires more computation time than \texttt{BARON}. However, given that we have implemented our algorithms in Python, it remains to be seen whether efficiency gains that would be achieved by utilizing a compiled language can provide for a more straightforward comparison.

\begin{table}[h]
 \caption{\label{tab:rootNodeStudy} Comparison of QCQP relaxation with \texttt{BARON} root node.}
 \centering 
 \begin{adjustbox}{width=0.75\textwidth}
 \begin{tabular}{c c c c c }
 \toprule 
  \multirow{2}{*}{Dimension} & \multirow{2}{*}{\# Problems} & \multirowcell{2}{\texttt{BARON} Root \\ Node Better} & \multirowcell{2}{QCQP Bound \\ Better} & \multirowcell{2}{Avg. Gap \\ Reduction by QCQP} \\
  & & & & \\
  \midrule
  1D & 6 & 1 & 5 & 78.8\%$^\ast$\\
  2D & 6 & 0 & 6 & 92.1\% \\
  3D & 6 & 0 & 6 & 94.4\% \\
  4D & 6 & 0 & 6 & 94.5\% \\
  \bottomrule
  \multicolumn{5}{l}{\small $^\ast$Reduction for the 5 problems where the QCQP provided a better bound.  }
 \end{tabular}
 \end{adjustbox}
 \end{table}

\section{Conclusions} \label{sec:conclusions}
In this work, we presented a hierarchy of methodologies to construct convex quadratic underestimators for non-convex d.c. functions.
Focusing on d.c. functions extracted from optimization benchmark libraries, we generated quadratic underestimators that reduce the hypervolume between the function and a linear underestimator constructed at the same point of construction by 55.9\% for one-dimensional functions and, depending on the variant of our methodology used, by a range of 38.4\%-50.2\% for two-dimensional functions, on average. Furthermore, we demonstrated in our computational study that, with variants that include a shift in the quadratic form, we could generate valid underestimators for all points of construction where $\nabla^2f(\x) \succcurlyeq 0$, including points at which the first-order Taylor series approximation is not a valid underestimator.
We showcase the tightness of our quadratic underestimators by providing qualitative results on several example functions, highlighting also the improvements that can be achieved by utilizing more involved methods in our hierarchy. Finally, we show the quality of convex QCQP relaxations constructed using our quadratic underestimation methodology in a comparison with \texttt{BARON} for lower bounds computed at the root node of a set of systematically created d.c. optimization problems. Notably, our convex QCQP relaxation is able to produce superior lower bound than those computed by \texttt{BARON} at the root node in the vast majority of cases, closing the root node relaxation gap by 90\%, on average.

Future work could investigate generalizing the hierarchy even further to include variants that determine a matrix for the second-order term of the Taylor series approximation from the entire positive semi-definite cone, rather than from a restricted subset limited to matrices satisfying diagonal dominance properties. In our experience, maintaining a quadratic that is monotonic in its parameters is much more involved with a parameter modifying the first-order term of the quadratic. Thus, future work could explore ways to tailor the cutting plane algorithm to achieve the monotonicity property in the context of modifying the first-order term, or altogether removing the monotonicity requirement from the algorithm while preserving efficiency and convergence.

\section*{Acknowledgments}
We acknowledge financial support from Mitsubishi Electric Research Labs (MERL) through the Center for Advanced Process Decision-making (CAPD) at Carnegie Mellon University. William Strahl also gratefully acknowledges support from the R.R.~Rothfus Graduate Fellowship in Chemical Engineering and the Chevron Graduate Fellowship in Chemical Engineering.

\section*{Appendix}\label{sec:appendA}

\renewcommand{\thefigure}{A\arabic{figure}}
\renewcommand{\thetable}{A\arabic{table}}
\setcounter{figure}{0}
\setcounter{table}{0}
The material supplied in this appendix provides the explicit function and problem formulations used in our computational study as well as supplemental details for some of the results. More specifically, Table~\ref{tab:COCONUTDCFunctions} displays the precise functions used in Section~\ref{sec:hierarchy_results}, while Tables~\ref{tab:hierarchyDetailsValid} and~\ref{tab:hierarchyDetailsInvalid} provide additional details for the computational results of the same section. 

\begin{figure}[htb!]
	\centering
	\includegraphics[width=.6\textwidth]{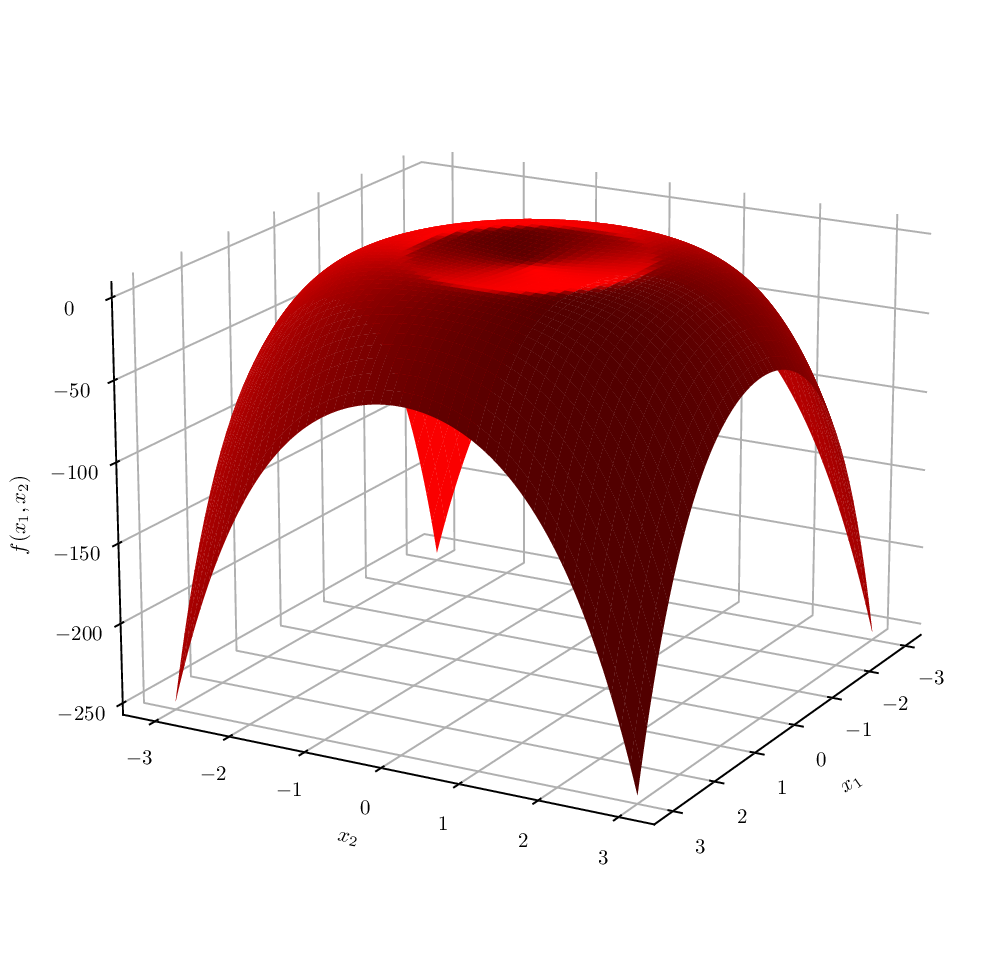}
	\caption{The pathological d.c. function from problem \textit{sisser}, where no locally convex point of construction admits a valid quadratic underestimator unless a method with shift is employed.}
	\label{fig:sisserPlot}
\end{figure}

Figure~\ref{fig:sisserPlot} provides a visual example of a pathological function discovered in the study for quadratic underestimator forms without the shift parameter.
Finally, for the computational study in Section \ref{sec:root_node_relax_comp}, where we compare root node relaxations of optimization problems, Tables \ref{tab:funcUsed}--\ref{tab:4D1rootNodeProblems} provide all the required information to construct the problems used in the study. All models are supplied in the Supplementary Material accompanying this paper in the form of \texttt{.nl} files, which also include the points used as points of construction for the computational results in Section~\ref{sec:root_node_comparison_results}.

\begin{table}[h!]
\caption{Details of functions from the COCONUT library used in our computational studies, where the equation names are as in the GAMS files available at \url{https://arnold-neumaier.at/glopt/coconut/Benchmark/Benchmark.html}.}
\label{tab:COCONUTDCFunctions}
\centering
 \adjustbox{width=\linewidth}{
 \begin{tabular}{ c c l l l l l l} 
  \toprule
  Dimension & Library & Problem & Name & Expression & Variable Bounds & Assigned Bounds\\
  \midrule
	1 & GlobalLib & ex4\_1\_6 & objcons & $f(x_1) = \left(27x_1^2 + x_1^6 +250\right) - 15x_1^4$ & (-5, 5) & -- \\
	1 & CUTE & zy2 & objcons & $f(x_1) = x_1^3 - 6x_1^2$ & (0, 8) & Upper \\
	1 & GlobalLib & ex4\_1\_9 & con1 & $f(x_1) = 8x_1^3 - \left(8x_1^2 + 2x_1^4\right)$ & (0, 3) & -- \\
	2 & CSTP & conform1 & con2, & $f(x_1, x_2) = \left[5x_1^2 + 5x_2^2 + \frac{3}{2}(x_1^2+x_2^2)^2\right] - \left[4(x_1 + x_2)^2 + \frac{3}{2}x_1^4 + \frac{3}{2}x_2^4\right]$ & (-3, 3), (-3, 3) & Both \\
	2 & GlobalLib & ex8\_1\_4 & objcons & $f(x_1, x_2) = \left[15x_1^2 + 9x_2^2 + x_1^6\right] - \left[3(x_1+x_2)^2 + 6.3x_1^4\right]$ & (-3, 3), (-3, 3) & Both \\
	2 & CUTE & camel6 & objcons & $f(x_1, x_2) = \left[\frac{7}{2}x_1^2 + 0.5(x_1+x_2)^2 + 4x_2^4 + \frac{1}{3}x_1^6\right] - \left(\frac{9}{2}x_2^2 + 2.1x_1^4\right)$ & (-3, 3), (-1.5, 1.5) & -- \\
	2 & CUTE & sisser & objcons & $f(x_1, x_2) = \left(4x_1^2 + 4x_2^2\right) - \left(x_1^2+x_2^2\right)^2$ & (-3, 3), (-3, 3) &  Both \\
	2 & CSTP & cyclo & con1 & $f(x_1, x_2) = \left[264.5x_1^2 + 79.5x_2^2 + 694.5(x_1 + x_2)^2 + 656.5\left(x_1^2+x_2^2\right)^2\right] - \left(656.5x_1^4 + 656.5x_2^4\right)$ & (-10, 10), (-10, 10) & -- \\
	2 & GlobalLib & ex4\_1\_5 & objcons & $f(x_1, x_2) = \left(\frac{3}{2}x_2^2 + \frac{5}{2}x_1^2 + \frac{1}{6}x_1^6\right) - \left[0.5(x_1+x_2)^2 + 1.05x_1^4\right]$ & (-5, 5), (-5, 5) & Upper for $x_1$ and lower for $x_2$ \\
	2 & CUTE & dipigri & objcons & $f(x_1, x_2) = \left(x_2^4 + 9x_1^2 + 2x_2^2\right) - 2(x_1+x_2)^2$ & (-3, 3), (-3, 3) & Both \\
   \bottomrule
 \end{tabular}}
 \end{table}

\begin{table}[h!]
	\caption{\label{tab:hierarchyDetailsValid} Computational details (averages $\pm$ standard deviations) for producing underestimators of bivariate functions using the 90 points of construction (in total, across all functions) that did not require shift for successful construction.}
	\centering
	\begin{tabular}{ c c c c c c } 
	 \toprule
	 Method & Iterations & Vertices & Metric & CPU (ms) & LP Solves \\
	 \midrule
 	  S & $44.0\pm31.2$ & $205.8\pm147.5$ & $0.384\pm0.289$ & $ 35\pm32$ & $0.0\pm0.0$\\ 
      D & $56.4\pm35.7$ & $265.0\pm168.8$ & $0.449\pm0.288$ & $ 93\pm65$ & $4.3\pm1.9$\\ 
      M & $61.0\pm38.2$ & $287.3\pm181.0$ & $0.464\pm0.293$ & $125\pm73$ & $4.4\pm1.9$\\ 
     SS & $44.0\pm31.2$ & $205.8\pm147.5$ & $0.384\pm0.289$ & $ 44\pm42$ & $0.0\pm0.0$\\ 
    UDS & $53.6\pm37.8$ & $251.5\pm179.8$ & $0.467\pm0.287$ & $ 98\pm68$ & $4.0\pm1.8$\\ 
     DS & $63.6\pm40.9$ & $299.3\pm194.6$ & $0.499\pm0.280$ & $114\pm75$ & $4.8\pm2.1$\\ 
     MS & $65.2\pm40.8$ & $307.2\pm194.6$ & $0.502\pm0.282$ & $133\pm68$ & $4.8\pm2.1$\\ 
	 \bottomrule
	\end{tabular}
\end{table}

\begin{table}[h!]
	\caption{\label{tab:hierarchyDetailsInvalid} Computational details (averages $\pm$ standard deviations) for producing underestimators of bivariate functions using the 85 points of construction (in total, across all functions) that required the shift for successful construction.}
	\centering
	\begin{tabular}{ c c c c c c } 
	 \toprule
	 Method & Iterations & Vertices & Metric & CPU (ms) & LP Solves \\
	 \midrule
	 SS & $33.1\pm24.6$ & $155.1\pm115.8$ & $0.000\pm0.000$ & $34\pm27$ & $0.0\pm0.0$\\ 
    UDS & $40.9\pm30.5$ & $192.6\pm146.9$ & $0.094\pm0.194$ & $75\pm49$ & $3.1\pm2.0$\\ 
     DS & $41.6\pm31.2$ & $195.8\pm149.9$ & $0.115\pm0.190$ & $73\pm56$ & $3.4\pm2.2$\\ 
     MS & $42.8\pm31.6$ & $201.1\pm151.7$ & $0.122\pm0.195$ & $96\pm47$ & $3.4\pm2.2$\\ 
	 \bottomrule
	\end{tabular}
\end{table}
		
\begin{table}[htb!]
\centering
	\caption{Functions utilized in root node study.}
	\label{tab:funcUsed}
	\begin{tabular}{l}
	\toprule
	Functions\\
	\midrule
	$f_1(x) := \left(20x^{10} + 4x^2\right) - 12x^4$\\
	$f_2(x) := \left(4x^6 + x^4\right) - 3x^2$\\
	$f_3(x) := \left(32x^6 + 8x^2\right) - 31x^4$\\
	$f_4(x) := \left(4x^2 + 8x^8\right) - \left(e^{2.35x} + e^{-2.35x} - 2.35\right)$\\
	$f_5(x) := \left(e^{4x}+e^{-4x}\right) - \left(6e^{2x} + 6e^{-2x} - 8.6\right)$\\
	$f_6(x) := \left(20x^{10} + 36x^6 +8x^2\right) - 38x^4$\\
	$\mathcal{L}_i(\x) := \frac{1}{2p}(\sum_{j=1}^p x_j)^4$\\
	\bottomrule
	\end{tabular}
\end{table}
	
\begin{landscape}
\begin{table}[htb!]
\centering
 \caption{1D d.c. optimization problems utilized in root node study.}
\label{tab:1D1rootNodeProblems}
\centering
\adjustbox{max width=\linewidth, max totalheight=\textheight}{
\begin{tabular}{@{} l l @{}}
 \toprule 
 Problem \# & Definition \\ 
 \hline
01 & {$\!\begin{array}{cl}
	\min\limits_{\mathbf{x} \in [-1,1]} & f_1(x_1)\\ 
	 \text{s.t.} 			  & 0.053 \left(0.909 x_{1}^{10} + 1\right)^{3} e^{1.0 x_{1}^{10}} \le 0.075\\ 
	  			 & f_2(x_1) \le -0.072\\ 
	\end{array}$}\\ 
\hline
02 & {$\!\begin{array}{cl}
	\min\limits_{\mathbf{x} \in [-1,1]} & f_6(x_1)\\ 
	 \text{s.t.} 			  & \frac{0.004}{\left(1 - \frac{0.091}{x_{1} + 1.1}\right)^{2.333}} \le 0.006\\ 
	  			 & f_2(x_1) \le -0.031\\ 
	  			 & f_5(x_1) \le 0.466\\ 
	\end{array}$}\\ 
\hline
03 & {$\!\begin{array}{cl}
	\min\limits_{\mathbf{x} \in [-1,1]} & f_4(x_1)\\ 
	 \text{s.t.} 			  & - 1.103 \left(1 - 0.334 e^{- x_{1}}\right)^{0.75} \le -0.540\\ 
	  			 & f_1(x_1) \le 0.320\\ 
	  			 & f_6(x_1) \le 0.381\\ 
	  			 & f_5(x_1) \le 0.579\\ 
	\end{array}$}\\ 
 \hline
04 & {$\!\begin{array}{cl}
	\min\limits_{\mathbf{x} \in [-1,1]} & f_4(x_1)\\ 
	 \text{s.t.} 			  & \frac{0.053 \left(- \frac{0.909 \log{\left(1.1 - x_{1} \right)}}{\log{\left(10 \right)}} + 1\right)^{3}}{\left(1.1 - x_{1}\right)^{\frac{1.0}{\log{\left(10 \right)}}}} \le 0.152\\ 
	  			 &- 1.0 \left(1 - 0.909 x_{1}^{10}\right)^{0.1} \le -0.998\\ 
	  			 & f_6(x_1) \le 0.422\\ 
	\end{array}$}\\ 
\hline
05 & {$\!\begin{array}{cl}
	\min\limits_{\mathbf{x} \in [-1,1]} & f_6(x_1)\\ 
	 \text{s.t.} 			  & \frac{0.202}{\left(- 0.048 \left(0.909 x_{1} + 1\right)^{3} e^{x_{1}} + 1\right)^{0.667}} \le 0.263\\ 
	  			 &- 0.616 \left(0.852 \left(0.909 x_{1} + 1\right)^{0.1} + 1\right)^{0.75} \le -0.946\\ 
	  			 & f_2(x_1) \le -0.059\\ 
	  			 & f_5(x_1) \le 0.271\\ 
	\end{array}$}\\ 
\hline
06 & {$\!\begin{array}{cl}
	\min\limits_{\mathbf{x} \in [-1,1]} & f_5(x_1)\\ 
	 \text{s.t.} 			  & \frac{0.1}{- \frac{0.074 e^{- x_{1}}}{\left(0.909 x_{1} + 1\right)^{0.667}} + 1.1} \le 0.115\\ 
	  			 &\frac{0.202}{\left(0.909 x_{1} + 1\right)^{0.667}} \le 0.247\\ 
	  			 & f_2(x_1) \le -0.018\\ 
	  			 & f_4(x_1) \le 0.305\\ 
	  			 & f_3(x_1) \le 0.616\\ 
	\end{array}$}\\ 
 \bottomrule 
\end{tabular}}
\end{table}

\begin{table}[htb!]
 \caption{2D d.c. optimization problems utilized in root node study.}
\label{tab:2D1rootNodeProblems}
\centering
\adjustbox{max width=\linewidth, max totalheight=\textheight}{
\begin{tabular}{ l l }
 \toprule 
 Problem \# & Definition \\ 
 \hline
07 & {$\!\begin{array}{cl}
	\min\limits_{\mathbf{x} \in [-1,1]^2} & f_2(x_1) + f_5(x_2) + \mathcal{L}_2(x_1, x_2)\\ 
	 \text{s.t.} 			  & - x_{1} + x_{2} \le 0.797\\ 
	 			 &\frac{1.887 e^{- 0.469 \left(0.909 x_{1} + 1\right)^{0.1} - 0.308 \left(0.909 x_{2} + 1\right)^{0.75}}}{\left(0.426 \left(0.909 x_{1} + 1\right)^{0.1} + 0.28 \left(0.909 x_{2} + 1\right)^{0.75} + 1\right)^{0.667}} \le 0.760\\ 
	  			 & f_4(x_1) + f_6(x_2) + \mathcal{L}_2(x_1, x_2) \le 0.730\\ 
	\end{array}$}\\ 
\hline
08 & {$\!\begin{array}{cl}
	\min\limits_{\mathbf{x} \in [-1,1]^2} & f_1(x_1) + f_5(x_2) + \mathcal{L}_2(x_1, x_2)\\ 
	 \text{s.t.} 			  & - x_{1} - x_{2} \le 0.688\\ 
	  			 &\frac{0.053 \left(\frac{0.034 e^{- x_{2}}}{\left(0.909 x_{2} + 1\right)^{0.667}} - \frac{0.455 \log{\left(1.1 - x_{1} \right)}}{\log{\left(10 \right)}} + 1\right)^{3} e^{\frac{0.037 e^{- x_{2}}}{\left(0.909 x_{2} + 1\right)^{0.667}}}}{\left(1.1 - x_{1}\right)^{\frac{0.5}{\log{\left(10 \right)}}}} \le 0.132\\ 
	  			 & f_3(x_1) + f_4(x_2) + \mathcal{L}_2(x_1, x_2) \le 0.963\\ 
	  			 & f_2(x_1) + f_6(x_2) + \mathcal{L}_2(x_1, x_2) \le 0.219\\ 
	\end{array}$}\\ 
\hline
09 & {$\!\begin{array}{cl}
	\min\limits_{\mathbf{x} \in [-1,1]^2} & f_1(x_1) + f_4(x_2) + \mathcal{L}_2(x_1, x_2)\\ 
	 \text{s.t.} 			  & x_{1} - x_{2} \le 0.727\\ 
	  			 &0.144 \left(0.024 \left(0.909 x_{2} + 1\right)^{3} e^{x_{2}} + 1 + \frac{0.045}{x_{1} + 1.1}\right)^{3} \le 0.307\\ 
	  			 & f_2(x_1) + f_3(x_2) + \mathcal{L}_2(x_1, x_2) \le 0.580\\ 
	  			 & f_1(x_1) + f_3(x_2) + \mathcal{L}_2(x_1, x_2) \le 0.888\\ 
	  			 & f_2(x_1) + f_4(x_2) + \mathcal{L}_2(x_1, x_2) \le 0.052\\ 
	\end{array}$}\\ 
 \hline
10 & {$\!\begin{array}{cl}
	\min\limits_{\mathbf{x} \in [-1,1]^2} & f_4(x_1) + f_5(x_2) + \mathcal{L}_2(x_1, x_2)\\ 
	 \text{s.t.} 			  & - x_{1} + x_{2} \le 0.741\\ 
	  			 &\frac{0.45 e^{- 0.308 \left(0.909 x_{1} + 1\right)^{0.75} + 0.184 e^{- x_{2}}}}{\left(0.28 \left(0.909 x_{1} + 1\right)^{0.75} + 1 - 0.167 e^{- x_{2}}\right)^{0.667}} \le 0.515\\ 
	  			 &\frac{0.074 e^{\frac{0.101}{\left(0.909 x_{1} + 1\right)^{0.667}} + 0.072 \left(0.909 x_{2} + 1\right)^{3}}}{\left(- \frac{0.092}{\left(0.909 x_{1} + 1\right)^{0.667}} - 0.065 \left(0.909 x_{2} + 1\right)^{3} + 1\right)^{0.667}} \le 0.133\\ 
	  			 & f_2(x_1) + f_3(x_2) + \mathcal{L}_2(x_1, x_2) \le 0.559\\ 
	\end{array}$}\\ 
\hline
11 & {$\!\begin{array}{cl}
	\min\limits_{\mathbf{x} \in [-1,1]^2} & f_2(x_1) + f_4(x_2) + \mathcal{L}_2(x_1, x_2)\\ 
	 \text{s.t.} 			  & x_{1} - x_{2} \le 0.719\\ 
	  			 &0.053 \left(0.065 \left(0.909 x_{1} + 1\right)^{3} + 0.167 e^{x_{2}} + 1\right)^{3} e^{0.072 \left(0.909 x_{1} + 1\right)^{3} + 0.184 e^{x_{2}}} \le 0.275\\ 
	  			 &- 0.967 \left(\frac{0.455 \left(- 1.478 x_{1} - 1.626\right) \log{\left(x_{1} + 1.1 \right)}}{\log{\left(10 \right)}} - 0.167 e^{x_{2}} + 1\right)^{0.75} \le -0.718\\ 
	  			 & f_1(x_1) + f_2(x_2) + \mathcal{L}_2(x_1, x_2) \le 0.364\\ 
	  			 & f_5(x_1) + f_5(x_2) + \mathcal{L}_2(x_1, x_2) \le 2.221\\ 
	\end{array}$}\\ 
\hline
12 & {$\!\begin{array}{cl}
	\min\limits_{\mathbf{x} \in [-1,1]^2} & f_4(x_1) + f_5(x_2) + \mathcal{L}_2(x_1, x_2)\\ 
	 \text{s.t.} 			  & x_{1} - x_{2} \le 0.810\\ 
	  			 &0.144 \left(0.455 x_{1} + 0.455 x_{2} + 1\right)^{3} \le 0.363\\ 
	  			 &- 1.017 \left(- \frac{0.002}{\left(0.909 x_{1} + 1\right)^{2.333}} + 1 - \frac{0.045}{x_{2} + 1.1}\right)^{0.75} \le -0.942\\ 
	  			 & f_1(x_1) + f_2(x_2) + \mathcal{L}_2(x_1, x_2) \le 0.234\\ 
	  			 & f_2(x_1) + f_3(x_2) + \mathcal{L}_2(x_1, x_2) \le 0.474\\ 
	  			 & f_1(x_1) + f_3(x_2) + \mathcal{L}_2(x_1, x_2) \le 0.768\\ 
	\end{array}$}\\ 
 \bottomrule 
\end{tabular}}
\end{table}

\begin{table}[htb!]
 \caption{3D d.c. optimization problems utilized in root node study.}
\label{tab:3D1rootNodeProblems}
\centering
\adjustbox{max width=\linewidth, max totalheight=\textheight}{
\begin{tabular}{ l l }
 \toprule 
 Problem \# & Definition \\ 
 \hline
13 & {$\!\begin{array}{cl}
	\min\limits_{\mathbf{x} \in [-1,1]^3} & f_2(x_1) + f_4(x_2) + f_5(x_3) + \mathcal{L}_3(x_1, x_2, x_3)\\ 
	 \text{s.t.} 			  & x_{1} - x_{2} + x_{3} \le 0.938\\ 
	  			 &\frac{0.679}{- 0.333 x_{3}^{4} + 0.312 \left(0.909 x_{1} + 1\right)^{0.1} - 0.048 \left(0.909 x_{2} + 1\right)^{3} + 1.1} \le 0.605\\ 
	  			 & f_4(x_1) + f_4(x_2) + f_4(x_3) + \mathcal{L}_3(x_1, x_2, x_3) \le 0.750\\ 
	\end{array}$}\\ 
\hline
14 & {$\!\begin{array}{cl}
	\min\limits_{\mathbf{x} \in [-1,1]^3} & f_1(x_1) + f_2(x_2) + f_6(x_3) + \mathcal{L}_3(x_1, x_2, x_3)\\ 
	 \text{s.t.} 			  & x_{1} + x_{2} + x_{3} \le 0.891\\ 
	  			 &\frac{0.074 e^{\frac{0.067}{\left(0.909 x_{1} + 1\right)^{0.667}} + 0.123 e^{x_{3}}}}{\left(1.1 - x_{2}\right)^{\frac{0.333}{\log{\left(10 \right)}}} \left(- \frac{0.061}{\left(0.909 x_{1} + 1\right)^{0.667}} - 0.111 e^{x_{3}} + \frac{0.303 \log{\left(1.1 - x_{2} \right)}}{\log{\left(10 \right)}} + 1\right)^{0.667}} \le 0.130\\ 
	  			 & f_1(x_1) + f_4(x_2) + f_4(x_3) + \mathcal{L}_3(x_1, x_2, x_3) \le 1.277\\ 
	  			 & f_4(x_1) + f_5(x_2) + f_6(x_3) + \mathcal{L}_3(x_1, x_2, x_3) \le 2.414\\ 
	\end{array}$}\\ 
\hline
15 & {$\!\begin{array}{cl}
	\min\limits_{\mathbf{x} \in [-1,1]^3} & f_1(x_1) + f_1(x_2) + f_6(x_3) + \mathcal{L}_3(x_1, x_2, x_3)\\ 
	 \text{s.t.} 			  & x_{1} - x_{2} + x_{3} \le 0.873\\ 
	  			 &\frac{0.004}{\left(- 0.044 \left(0.909 x_{1} + 1\right)^{3} - \frac{0.023 e^{- x_{3}}}{\left(0.909 x_{3} + 1\right)^{0.667}} + \frac{0.303 \log{\left(1.1 - x_{2} \right)}}{\log{\left(10 \right)}} + 1\right)^{2.333}} \le 0.008\\ 
	  			 & f_2(x_1) + f_3(x_2) + f_4(x_3) + \mathcal{L}_3(x_1, x_2, x_3) \le 1.158\\ 
	  			 & f_4(x_1) + f_4(x_2) + f_4(x_3) + \mathcal{L}_3(x_1, x_2, x_3) \le 0.604\\ 
	  			 & f_3(x_1) + f_4(x_2) + f_4(x_3) + \mathcal{L}_3(x_1, x_2, x_3) \le 0.887\\ 
	\end{array}$}\\ 
 \hline
16 & {$\!\begin{array}{cl}
	\min\limits_{\mathbf{x} \in [-1,1]^3} & f_1(x_1) + f_3(x_2) + f_4(x_3) + \mathcal{L}_3(x_1, x_2, x_3)\\ 
	 \text{s.t.} 			  & - x_{1} - x_{2} + x_{3} \le 0.879\\ 
	  			 &- 0.975 \left(0.187 \left(0.909 x_{1} + 1\right)^{0.75} - 0.044 \left(0.909 x_{2} + 1\right)^{3} + 1 - \frac{0.03}{x_{3} + 1.1}\right)^{0.1} \le -0.969\\ 
	  			 &0.053 \left(0.044 \left(0.909 x_{1} + 1\right)^{3} + \frac{0.023 e^{- x_{3}}}{\left(0.909 x_{3} + 1\right)^{0.667}} + 1 + \frac{0.03}{x_{2} + 1.1}\right)^{3} e^{0.048 \left(0.909 x_{1} + 1\right)^{3} + \frac{0.025 e^{- x_{3}}}{\left(0.909 x_{3} + 1\right)^{0.667}} + \frac{0.033}{x_{2} + 1.1}} \le 0.158\\ 
	  			 & f_5(x_1) + f_5(x_2) + f_5(x_3) + \mathcal{L}_3(x_1, x_2, x_3) \le 4.602\\ 
	\end{array}$}\\ 
\hline
17 & {$\!\begin{array}{cl}
	\min\limits_{\mathbf{x} \in [-1,1]^3} & f_1(x_1) + f_1(x_2) + f_3(x_3) + \mathcal{L}_3(x_1, x_2, x_3)\\ 
	 \text{s.t.} 			  & - x_{1} - x_{2} - x_{3} \le 0.867\\ 
	  			 &0.053 \left(\frac{0.061}{\left(0.909 x_{1} + 1\right)^{0.667}} + 1 + 0.111 e^{- x_{2}} + \frac{0.03}{x_{3} + 1.1}\right)^{3} e^{\frac{0.067}{\left(0.909 x_{1} + 1\right)^{0.667}} + 0.123 e^{- x_{2}} + \frac{0.033}{x_{3} + 1.1}} \le 0.160\\ 
	  			 &0.867 \left(- 0.284 \left(0.909 x_{1} + 1\right)^{0.1} - 0.187 \left(0.909 x_{2} + 1\right)^{0.75} + 0.111 e^{x_{3}} + 1\right)^{3} \le 0.360\\ 
	  			 & f_1(x_1) + f_6(x_2) + f_6(x_3) + \mathcal{L}_3(x_1, x_2, x_3) \le 3.062\\ 
	  			 & f_1(x_1) + f_2(x_2) + f_5(x_3) + \mathcal{L}_3(x_1, x_2, x_3) \le 0.837\\ 
	\end{array}$}\\ 
\hline
18 & {$\!\begin{array}{cl}
	\min\limits_{\mathbf{x} \in [-1,1]^3} & f_2(x_1) + f_2(x_2) + f_3(x_3) + \mathcal{L}_3(x_1, x_2, x_3)\\ 
	 \text{s.t.} 			  & - x_{1} - x_{2} + x_{3} \le 0.855\\ 
	  			 &0.531 e^{- 0.205 \left(0.909 x_{1} + 1\right)^{0.75} + \frac{0.001}{\left(0.909 x_{2} + 1\right)^{2.333}} + \frac{0.033}{x_{3} + 1.1}} \le 0.501\\ 
	  			 &\frac{0.917 e^{- 0.312 \left(0.909 x_{1} + 1\right)^{0.1} - 0.205 \left(0.909 x_{2} + 1\right)^{0.75} - \frac{0.333 \left(- 1.478 x_{3} - 1.626\right) \log{\left(x_{3} + 1.1 \right)}}{\log{\left(10 \right)}}}}{\left(0.284 \left(0.909 x_{1} + 1\right)^{0.1} + 0.187 \left(0.909 x_{2} + 1\right)^{0.75} + \frac{0.303 \left(- 1.478 x_{3} - 1.626\right) \log{\left(x_{3} + 1.1 \right)}}{\log{\left(10 \right)}} + 1\right)^{0.667}} \le 0.511\\ 
	  			 & f_3(x_1) + f_5(x_2) + f_6(x_3) + \mathcal{L}_3(x_1, x_2, x_3) \le 5.649\\ 
	  			 & f_5(x_1) + f_5(x_2) + f_5(x_3) + \mathcal{L}_3(x_1, x_2, x_3) \le 3.941\\ 
	  			 & f_3(x_1) + f_6(x_2) + f_6(x_3) + \mathcal{L}_3(x_1, x_2, x_3) \le 1.192\\ 
	\end{array}$}\\
\bottomrule 
\end{tabular}}
\end{table}

\begin{table}[htb!]
 \caption{4D d.c. optimization problems utilized in root node study.}
\label{tab:4D1rootNodeProblems}
\centering
\adjustbox{max width=\linewidth, max totalheight=\textheight}{
\begin{tabular}{ l l }
 \toprule 
 Problem \# & Definition \\ 
 \hline
19 & {$\!\begin{array}{cl}
	\min\limits_{\mathbf{x} \in [-1,1]^4} & f_1(x_1) + f_1(x_2) + f_1(x_3) + f_4(x_4) + \mathcal{L}_4(x_1, x_2, x_3, x_4)\\ 
	 \text{s.t.} 			  & x_{1} - x_{2} + x_{3} - x_{4} \le 1.016\\ 
	  			 &- \frac{1.0 \log{\left(- \frac{0.001}{\left(0.909 x_{1} + 1\right)^{2.333}} - 0.013 \left(0.909 x_{4} + 1\right)^{3} e^{x_{4}} + 1.1 - 0.092 e^{- x_{2}} - \frac{0.025}{x_{3} + 1.1} \right)}}{\log{\left(10 \right)}} \le 0.089\\ 
	  			 & f_2(x_1) + f_2(x_2) + f_3(x_3) + f_5(x_4) + \mathcal{L}_4(x_1, x_2, x_3, x_4) \le 3.270\\ 
	\end{array}$}\\ 
\hline
20 & {$\!\begin{array}{cl}
	\min\limits_{\mathbf{x} \in [-1,1]^4} & f_1(x_1) + f_5(x_2) + f_6(x_3) + f_6(x_4) + \mathcal{L}_4(x_1, x_2, x_3, x_4)\\ 
	 \text{s.t.} 			  & x_{1} - x_{2} + x_{3} - x_{4} \le 0.906\\ 
	  			 &- 0.976 \left(0.14 \left(0.909 x_{1} + 1\right)^{0.75} + \frac{0.227 \left(- 1.478 x_{2} - 1.626\right) \log{\left(x_{2} + 1.1 \right)}}{\log{\left(10 \right)}} - \frac{0.001}{\left(0.909 x_{3} + 1\right)^{2.333}} - \frac{0.017 e^{- x_{4}}}{\left(0.909 x_{4} + 1\right)^{0.667}} + 1\right)^{0.1} \le -0.969\\ 
	  			 & f_1(x_1) + f_1(x_2) + f_2(x_3) + f_2(x_4) + \mathcal{L}_4(x_1, x_2, x_3, x_4) \le 1.404\\ 
	  			 & f_2(x_1) + f_2(x_2) + f_2(x_3) + f_3(x_4) + \mathcal{L}_4(x_1, x_2, x_3, x_4) \le 0.157\\ 
	\end{array}$}\\ 
\hline
21 & {$\!\begin{array}{cl}
	\min\limits_{\mathbf{x} \in [-1,1]^4} & f_1(x_1) + f_4(x_2) + f_5(x_3) + f_5(x_4) + \mathcal{L}_4(x_1, x_2, x_3, x_4)\\ 
	 \text{s.t.} 			  & x_{1} + x_{2} - x_{3} + x_{4} \le 0.988\\ 
	  			 &0.053 \left(0.033 \left(0.909 x_{1} + 1\right)^{3} + \frac{0.227 \cdot \left(1.478 x_{2} + 1.626\right) \log{\left(x_{2} + 1.1 \right)}}{\log{\left(10 \right)}} + 1 + 0.084 e^{- x_{3}} + \frac{0.023}{x_{4} + 1.1}\right)^{3} e^{0.036 \left(0.909 x_{1} + 1\right)^{3} + \frac{0.25 \cdot \left(1.478 x_{2} + 1.626\right) \log{\left(x_{2} + 1.1 \right)}}{\log{\left(10 \right)}} + 0.092 e^{- x_{3}} + \frac{0.025}{x_{4} + 1.1}} \le 0.151\\ 
	  			 & f_3(x_1) + f_4(x_2) + f_4(x_3) + f_4(x_4) + \mathcal{L}_4(x_1, x_2, x_3, x_4) \le 2.414\\ 
	  			 & f_2(x_1) + f_3(x_2) + f_3(x_3) + f_6(x_4) + \mathcal{L}_4(x_1, x_2, x_3, x_4) \le 2.160\\ 
	  			 & f_2(x_1) + f_2(x_2) + f_2(x_3) + f_5(x_4) + \mathcal{L}_4(x_1, x_2, x_3, x_4) \le -0.047\\ 
	\end{array}$}\\ 
\hline
22 & {$\!\begin{array}{cl}
	\min\limits_{\mathbf{x} \in [-1,1]^4} & f_1(x_1) + f_3(x_2) + f_3(x_3) + f_3(x_4) + \mathcal{L}_4(x_1, x_2, x_3, x_4)\\ 
	 \text{s.t.} 			  & - x_{1} + x_{2} - x_{3} - x_{4} \le 0.938\\ 
	  			 &- 0.969 \left(- 0.227 x_{3}^{4} - \frac{0.046}{\left(0.909 x_{1} + 1\right)^{0.667}} - 0.033 \left(0.909 x_{2} + 1\right)^{3} + \frac{0.227 \log{\left(1.1 - x_{4} \right)}}{\log{\left(10 \right)}} + 1\right)^{0.75} \le -0.770\\ 
	  			 &0.219 \left(0.227 x_{3}^{10} - 0.14 \left(0.909 x_{1} + 1\right)^{0.75} + 0.012 \left(0.909 x_{4} + 1\right)^{3} e^{x_{4}} + 1 + \frac{0.023}{x_{2} + 1.1}\right)^{3} \le 0.228\\ 
	  			 & f_2(x_1) + f_2(x_2) + f_2(x_3) + f_6(x_4) + \mathcal{L}_4(x_1, x_2, x_3, x_4) \le 0.035\\ 
	\end{array}$}\\ 
\hline
23 & {$\!\begin{array}{cl}
	\min\limits_{\mathbf{x} \in [-1,1]^4} & f_2(x_1) + f_2(x_2) + f_4(x_3) + f_5(x_4) + \mathcal{L}_4(x_1, x_2, x_3, x_4)\\ 
	 \text{s.t.} 			  & - x_{1} - x_{2} - x_{3} + x_{4} \le 1.062\\ 
	  			 &\frac{0.488}{\left(- 0.227 x_{2}^{4} - 0.227 x_{4}^{10} + 0.14 \left(0.909 x_{1} + 1\right)^{0.75} + 1 - 0.084 e^{- x_{3}}\right)^{0.667}} \le 0.520\\ 
	  			 &\frac{0.074 e^{\frac{0.051}{\left(0.909 x_{1} + 1\right)^{0.667}} + 0.013 \left(0.909 x_{4} + 1\right)^{3} e^{x_{4}} + 0.092 e^{x_{2}} + 0.092 e^{- x_{3}}}}{\left(- \frac{0.046}{\left(0.909 x_{1} + 1\right)^{0.667}} - 0.012 \left(0.909 x_{4} + 1\right)^{3} e^{x_{4}} - 0.084 e^{x_{2}} + 1 - 0.084 e^{- x_{3}}\right)^{0.667}} \le 0.140\\ 
	  			 & f_1(x_1) + f_4(x_2) + f_5(x_3) + f_5(x_4) + \mathcal{L}_4(x_1, x_2, x_3, x_4) \le 8.879\\ 
	  			 & f_3(x_1) + f_3(x_2) + f_3(x_3) + f_4(x_4) + \mathcal{L}_4(x_1, x_2, x_3, x_4) \le 1.999\\ 
	\end{array}$}\\ 
\hline
24 & {$\!\begin{array}{cl}
	\min\limits_{\mathbf{x} \in [-1,1]^4} & f_1(x_1) + f_3(x_2) + f_5(x_3) + f_5(x_4) + \mathcal{L}_4(x_1, x_2, x_3, x_4)\\ 
	 \text{s.t.} 			  & x_{1} + x_{2} + x_{3} - x_{4} \le 1.016\\ 
	  			 &- \frac{1.0 \log{\left(- 0.25 x_{4}^{10} - \frac{0.051}{\left(0.909 x_{1} + 1\right)^{0.667}} - 0.036 \left(0.909 x_{2} + 1\right)^{3} + 1.1 - \frac{0.025}{x_{3} + 1.1} \right)}}{\log{\left(10 \right)}} \le 0.081\\ 
	  			 &- 0.998 \left(- 0.227 x_{4}^{10} + \frac{0.227 \left(- 1.478 x_{1} - 1.626\right) \log{\left(x_{1} + 1.1 \right)}}{\log{\left(10 \right)}} - \frac{0.001}{\left(0.909 x_{2} + 1\right)^{2.333}} + 1 - 0.084 e^{- x_{3}}\right)^{0.1} \le -0.970\\ 
	  			 & f_1(x_1) + f_1(x_2) + f_1(x_3) + f_3(x_4) + \mathcal{L}_4(x_1, x_2, x_3, x_4) \le 3.115\\ 
	  			 & f_1(x_1) + f_3(x_2) + f_5(x_3) + f_6(x_4) + \mathcal{L}_4(x_1, x_2, x_3, x_4) \le 2.389\\ 
	  			 & f_3(x_1) + f_6(x_2) + f_6(x_3) + f_6(x_4) + \mathcal{L}_4(x_1, x_2, x_3, x_4) \le 1.444\\ 
	\end{array}$}\\ 
 \bottomrule 
\end{tabular}}
\end{table}
\end{landscape}

\clearpage
\bibliographystyle{plain}
\bibliography{StrahlEtal_2_DC}
\end{document}